\documentclass[twocolumn]{autart}    

\usepackage{graphicx}     
\usepackage{amsfonts,latexsym,amssymb,amsmath}
\usepackage{url,stackengine}
\usepackage{color}

\catcode`\@=11
\def\downparenfill{$\m@th\braceld\leaders\vrule\hfill\bracerd$}
\def\downparenfill{$\m@th\braceld\leaders\vrule\hfill\bracerd$}
\def\overparen#1{\mskip 2mu\mathop{\vbox{\ialign{##\crcr\crcr \noalign{\kern0.4ex}
\downparenfill\crcr\noalign{\kern0.4ex\nointerlineskip}
$\hfil\displaystyle{#1}\hfil$\crcr}}}\limits\mskip 2mu} 
\catcode`\@=12

\newtheorem{assumption}{Assumption}
\newtheorem{lemma}{Lemma}
\newtheorem{theorem}{Theorem}
\newtheorem{definition}{Definition}

\newtheorem{remark}{Remark}
\newtheorem{example}{Example}
\newtheorem{proposition}{Proposition}

\newcommand{\NN}{\mathbb{N}}

\newcommand{\RR}{\mathbb{R}}

\newcommand\BB{\mathbb{B}} 
\def\X{\mathcal{X}}  
\def\T{\mathcal{T}}
\def\Z{\mathcal{Z}}

\def\O{\mathcal{O}}

\def\C{\mathcal{C}}  
\def\A{\mathcal{A}}  
\def\B{\mathcal{B}}

\def\H{\mathcal{H}}
\def\G{\mathcal{G}}
\def\P{\mathcal{P}}

\DeclareMathOperator{\card}{card}

\DeclareMathOperator{\diag}{diag}
\DeclareMathOperator*{\argmin}{argmin}

\DeclareMathOperator{\cl}{cl}

\def\ind{\mathcal{I}}
\def\indist{\,\stackunder{$\thicksim$}{}\,}
\def\nind{p}
\def\Tinv{T^{\rm inv}}

\begin{document}

\begin{frontmatter}

\title{Reconstructing  Indistinguishable Solutions \\ Via Set-Valued  KKL Observer}

\author[First]{Pauline Bernard}\ead{pauline.bernard@minesparis.psl.eu}, 
\author[Second]{Mohamed Maghenem}\ead{mohamed.maghenem@gipsa-lab.fr} 

\address[First]{Centre Automatique et Systèmes, Mines Paris, Université PSL, 60 bd Saint-Michel, 75272, Paris, France}
       
\address[Second]{University of Grenoble Alpes, CNRS, Gipsa-lab,  Grenoble INP,  France.}
          

\begin{abstract}                          
KKL observer design consists in finding a smooth change of coordinates transforming the system dynamics into a linear filter of the output. The state of the original system is then reconstructed by implementing this filter from any initial condition and left-inverting the transformation,  under a \textit{backward-distinguishability} property. In this paper, we consider the case where the latter assumption does not hold, namely when distinct solutions may generate the same output, and thus be indistinguishable. The KKL transformation is no longer injective and its ``left-inverse'' is thus allowed to be set-valued, yielding a set-valued KKL observer. Assuming the transformation is full-rank and its preimage has constant cardinality, we show the existence of a globally defined set-valued  left-inverse that is Lipschitz in the Hausdorff sense. Leveraging on recent results linking this left-inverse with the \textit{backward-indistinguishable sets}, we show that the set-valued KKL observer converges in the Hausdorff sense to the backward-indistinguishable set of the system solution. When, additionally, a given output is generated by a specific number of solutions not converging to each other, we show that the designed observer asymptotically reconstructs each of those solutions. Finally, the different assumptions are discussed and illustrated via examples.  
\end{abstract}

\end{frontmatter}

\section{Introduction}

Consider a dynamical system  of the form 
\begin{equation}
\label{sys}
 \dot x = f(x),  \qquad 	 y = h(x), 
\end{equation}
where  $x \in \RR^{n_x}$ is the state,   $y\in \RR^{n_y}$ is the output, $f: \RR^{n_x}\to \RR^{n_x}$ is assumed to be locally Lipschitz and $h:\RR^{n_x} \to \RR^{n_y}$ continuous.  
We assume that the trajectories of interest remain in a compact subset $\X \subset \RR^{n_x}$. A typical problem in many engineering applications is to estimate online the current state $x(t)$ of \eqref{sys} based on the knowledge of the output $y$ on the interval $[0,t]$, in the sense that, the error between the estimate, denoted by $\hat{x}(t)$, and $x(t)$ asymptotically converges to zero. To do so, we usually design an \textit{observer}; namely,  a dynamical system of the form 
\begin{equation} \label{obs}
\dot z = F(z,y),  \quad \quad	 \hat x = \T(z,y), 
\end{equation}
fed with the known output $y$, which provides, as output, an estimate $\hat{x}$ such that $\lim_{t\to \infty} (\hat{x}(t)-x(t)) = 0$ (see  \cite{bernardObserverDesignContinuoustime2022}). Since the observer cannot distinguish solutions producing a same output, the possibility of achieving asymptotic estimation implicitly requires that  \textit{indistinguishable} solutions, i.e., solutions producing a same output,  at least asymptotically converge to each other. In other words, the system must 
be \textit{detectable}. In general, nonlinear observers are designed under stronger \textit{observability} conditions
saying that there is no indistinguishable solutions, or, in other words, that the output signal  carries enough information to determine the state uniquely.  

\subsection{Background on KKL-Observer design}

One possible route to design \eqref{obs} is the so called \textit{nonlinear Luenberger} (or \textit{Kazantzis-Kravaris-Luenberger (KKL)}) approach,  initially introduced in  
  \cite{luenbergerObservingStateLinear1964} for single-output linear systems.  
The idea is to look for a continuous and injective change of coordinates $T\in \RR^{n_x\times n_x}$ such that $z := T(x)$ is governed by
\begin{equation}  \label{z_system}
\dot z = A z + B y,
\end{equation}
for a pair $(A,B)\in \RR^{n_x\times n_x} \times \RR^{n_x}$ to be chosen with $A$ Hurwitz.
In particular,  when $T$ is smooth, it must verify
\begin{equation} \label{eq_PDE}
\frac{\partial T}{\partial x}(x) f(x) = A T(x) + Bh(x) \qquad \forall x \in \X.
\end{equation}
Since $A$ is Hurwitz,  any solution to \eqref{z_system} exponentially estimates $T(x)$, and an estimate of $x$ can be recovered through left-inversion of $T$.

Originally, when $n_y=1$, D. Luenberger showed in \cite{luenbergerObservingStateLinear1964} that such a map $T$ always exists for linear observable systems as long as $(A,B)$ is picked controllable, of dimension $n_z=n_x$, and $A$ does not share any eigenvalue with the system's  dynamics. In the context of nonlinear systems, the existence of the map $T$ was first established around an equilibrium point in \cite{Shoshitaishvili_TSA_90}, \cite{kazantzisNonlinearObserverDesign1998} and \cite{krener2001nonlinear}.  Then, the localness was relaxed in \cite{KreisselmeierEngel_TAC_03} using a strong observability assumption,  which,  unfortunately,  does not provide an indication on the necessary dimension $n_z$ of $z$. This problem is solved in \cite{andrieuExistenceKazantzisKravaris2006a} by proving the existence of the injective map $T$,  under a weak \textit{backward-distinguishability} condition, for $A$ complex diagonal of dimension $n_x + 1$,  with a generic choice of  $n_x+1$ \textit{distinct} complex eigenvalues.  The aforementioned result is generalized in \cite{BriAndBerSer},  for \textit{almost any} real controllable pair $(A,B)$ of dimension $n_z=2n_x+1$ with $A$ diagonalizable. 

The distinguishability property assumed in \cite{andrieuExistenceKazantzisKravaris2006a} and \cite{BriAndBerSer} requires that the backward solutions from any distinct states $x_a,x_b$ in $\X$ generate distinct past outputs $y_a, y_b$; namely, we say that $x_a$ and $x_b$ are \textit{backward-distinguishable}.  Said differently, any given state in $\X$ is uniquely determined by its past output. 
However, this assumption is not always verified in applications, where some systems may exhibit indistinguishable states. In which case, estimating the state is typically impossible. However,  developing strategies to recover one of the indistinguishable states, or to observe all the possible trajectories corresponding to a given output, is still of great interest as we explain next.    

\subsection{Context and Motivation }

When two solutions, not asymptotically tending to each other, generate a same output,  there is no hope to design an observer producing a single 
asymptotically-correct estimate. However, one could imagine to have an algorithm producing  a set of estimates that  \textit{converge} asymptotically  --- for a certain distance to be defined --- to the set of solutions generating that same output, or producing one estimate converging asymptotically to one of the several possible solutions.   

The interest in such observers is motivated by classes of nonlinear systems producing \textit{finite} numbers of indistinguishable solutions.  This is illustrated in \cite{MorMujEsp,VerLorForMenZar} in the context of induction motors,  and in \cite{BerPra21} in the context of permanent magnet synchronous motors (PMSM)s.   Indeed,  when the number of possible indistinguishable solutions corresponding to a given output signal is finite,   one could hope to design an observer whose output $\hat{x}$ reconstructs one of the possible indistinguishable trajectories.   This is done in \cite{MorMujEsp,MorBes},   through sliding mode tools,  for systems that  can be written in an ``observable-like'' form. 
Instead,  in \cite{BerPra21},  the KKL design is used on a particular application featuring a PMSM with unknown resistance.  Indeed,  it is shown that indistinguishable solutions exist,  but there are always less than six,  and that there exists a map $T:\RR^{n_x} \to \RR^{n_z}$ transforming the dynamics into \eqref{z_system} and  whose inversion enables to reconstruct all the possible states.  The preliminary work in \cite{PaulineNolcos} attempted to generalize this idea by showing that a KKL observer is able to extract from a given output signal $t \mapsto y(t)$ all the possible information about the corresponding (indistinguishable) solutions. More precisely, it is proved that, for an  appropriate choice of the pair $(A,B)$,  there exists a continuously differentiable transformation $T : \mathcal{X} \rightarrow \mathbb{R}^{n_z}$ that transforms \eqref{sys} into the form of  \eqref{z_system} and that the map $T$ \textit{characterizes the distinguishable states}, in the sense that, its preimage gives exactly the set of indistinguishable states. This suggests that, when the preimage map is continuous (in the sense of set-valued maps), then it might be possible to online reconstruct, from  a solution $z$ to \eqref{z_system} subject to an output $y$, all the possible indistinguishable solutions generating $y$. Although this was confirmed on a particular simulation example, unless further assumptions are made, the preimage map is only defined on the image set $T(\X)$ and only upper semicontinuous in general, which prevents from stating a general convergence result.

\subsection{Contribution}

In this paper, we push forward the theory of set-valued KKL observers when the backward-distinguishability assumption is not satisfied. More specifically, we consider a general nonlinear system \eqref{sys} having a finite and constant number of indistinguishable states on $\X$ (resp. solutions). A smooth map $T$ transforms it into the form of \eqref{z_system}. As a result, the (set-valued) preimage map of $T$, denoted by $T^- : T(\X) \rightrightarrows \X$, allows us to generate what we call a set-valued KKL observer for \eqref{sys}.   Our goal is to propose sufficient conditions to guarantee that the designed observer asymptotically reconstructs the  possible solutions generating a same output. The paper's contributions can be listed as follows:

1- We show that the map $T^-$ is Lipschitz continuous provided that the jacobian of $T$ is full rank on $\X$  and $T^-$ has a constant cardinality on $T(\X)$.

2- Since  the solutions to \eqref{z_system} are not guaranteed to remain in $T(\X)$, we prove the existence of a Lipschitz continuous set-valued map $\Tinv : \mathbb{R}^{n_z} \rightrightarrows \mathbb{R}^{n_x}$ that is an extension of $T^-$ to $\mathbb{R}^{n_z}$. 

3- As a consequence of the latter two items, we show that, for any $x$ solution to \eqref{sys} generating the output $y$ and for any $z$ solution to \eqref{z_system} subject to $y$, 
the Hausdorff distance between
$\Tinv(z(t))$ and  the backward-indistinguishable set of $x(t)$ converges to zero.  

4- We provide further assumptions  under which,   any continuous selection $t \mapsto \hat{x}(t) \in \Tinv (z(t))$   converges asymptotically to a solution to \eqref{sys} generating the output $y$.  In particular, we consider the case where the cardinality of $T^-$ equals the number of solutions generating the output signal $y$ and eventually remaining in $\X$.

5-  Finally, we establish some connections between the different assumptions used in this paper. In particular, we investigate the link between the existence of a finite number of indistinguishable trajectories, and the constant and finite cardinality of the indistinguishable sets. We also relate the rank of $T$ with the rank of a \textit{differential observability map}, that can be more easily computed.

The rest of this paper is organized as follows. Preliminaries on indistinguishability and  continuity notions for set-valued maps are given in Section \ref{Sec.II} and the problem is stated in Section \ref{Sec.III}. Lipschitz continuity and Lipschitz extension of the preimage map $T^-$ are then analysed in Section \ref{Sec.IV}, while convergence properties of the proposed set-valued KKL observer are presented in Section \ref{sec_conv_KKL}. The link between the different assumptions is finally investigated in Section \ref{Sec.VI}. Examples are presented all along the paper to illustrate our results.

\textbf{Notations.} 
For $x \in \mathbb{R}^{n_x}$,  
$|x|$ denotes the Euclidean norm of $x$. For a set $K \subset \mathbb{R}^{n_x}$, we use $\mbox{int}(K)$ to denote its interior, $\partial K$ its boundary,  $\card (K)$  the number of elements in the set $K$, and, for a point $x\in\RR^{n_x}$, $d(x,K)$ denotes the distance from $x$ to the set $K$. 
For $O \subset \mathbb{R}^{n_x}$,  $K \backslash O$ denotes the subset of elements in $K$ that are not in $O$. The set $K$ is simply connected if any loop in $K$ can be continuously contracted to a point.  The \textit{Dubovitsky-Miliutin} cone of $K$ at $x$ is given by
\begin{align} \label{eq.cone2}
D_K(x) := & \{v \in \mathbb{R}^n : \exists \epsilon >0: \nonumber \\ &
x + \delta (v + w) \in K~\forall \delta \in (0,\epsilon],~\forall w \in \epsilon \mathbb{B} \},
\end{align}
and the \textit{contingent} cone of $K$ at $x$ is given by
\begin{align} \label{eq.toncon} 
T_K(x) := \left\{ v \in \mathbb{R}^n: \liminf_{h \rightarrow 0^+} \frac{|x + h v|_K}{h} = 0 \right\}.
\end{align}
For a differentiable map $x \mapsto T(x) \in \mathbb{R}^{n_z}$,  
$\frac{\partial T}{\partial x}$ denotes the Jacobian of $T$ with respect to $x$.  For $\varepsilon>0$, we use $\B(x,\varepsilon)$ to denote the open ball centered at $x$ whose radius is $\varepsilon$. By $F : \mathbb{R}^{n_z} \rightrightarrows \mathbb{R}^{n_x}$,  we denote a set-valued map associating to each element $z \in \mathbb{R}^{n_z}$ a subset $F(z) \subset \mathbb{R}^{n_x}$, and, for some $\Z \subset \mathbb{R}^{n_z}$, we let $F(\Z) := \{ \eta \in F(z) : z \in \Z \} = \cup_{z \in \Z}F(z)$. The set-valued map $F$ is said to be $p$-valued, for some $p \in \{1,2,... \}$, if for all $z \in \mathbb{R}^{n_z}$, $\card(F(z))=p$.
Furthermore, we denote $\A_p(\RR^{n_x})$  the space of unordered
$\nind$-tuples in $\RR^{n_x}$. Note that multiplicity is allowed in $\A_p(\RR^{n_x})$, i.e., a map $\tilde{F} : \mathbb{R}^{n_z} \rightarrow \A_p(\RR^{n_x})$ associates to each $z \in \mathbb{R}^{n_z}$, a tuple of $p$ points in $\RR^{n_x}$ that may not be distinct; in which case, we say that $\tilde{F}$ is an \textit{Almgren} $p$-valued map \cite{Alm}. Given $x_0\in \RR^n$, we denote $\omega(x_0)$ the $\omega$-limit set of $x_0$ by the dynamical system \eqref{sys}, defined by
$$ \omega(x_0) := \{ y  : \forall \bar{t}>0, ~ \forall \varepsilon > 0,~ \exists t > \bar{t} : |\phi(t) - y| < \varepsilon \}, $$
where $\phi$ is the unique maximal solution to \eqref{sys} starting from $x_0$.

\section{Preliminaries} \label{Sec.II}
Before stating the problem tackled in this paper, we recall some important related notions.

\subsection{Indistinguishability}

\begin{definition}[Backward indistinguishable points] \label{Def-dis-dis}
Two points  $x_a,x_b\in \RR^{n_x}$ are said to be  backward indistinguishable for \eqref{sys},   which we denote $x_a \indist x_b$, if the maximal backward solutions $\phi_a:(-\sigma_a,0]\to \RR^{n_x}$, $\phi_b:(-\sigma_b,0]\to \RR^{n_x}$ to \eqref{sys} initialized at $x_a,x_b$ respectively, i.e. with
$$
\phi_a(0) = x_a, \quad \quad \phi_b(0) = x_b,
$$
verify
$$
h(\phi_a(t)) = h(\phi_b(t)) \qquad \forall t \in (-\min\{\sigma_a,\sigma_b\}, 0] \ . 
$$
\end{definition}

In other words, two points $x_a$ and $x_b$ are backward indistinguishable for \eqref{sys} if they cannot be distinguished from the past output of \eqref{sys}. Indistinguishability may be seen as an equivalence relation whose classes of equivalence define the \textit{indistinguishable sets}.

\begin{definition}
[Indistinguishability set]
Given $\O\subseteq\RR^{n_x}$ and $x\in \RR^{n_x}$,  the backward indistinguishable set with respect to $\O$ for \eqref{sys} is given by
\begin{equation}
\label{eq_indist}
\ind_\O(x) := \{ x' \in \O, \: x'\indist x \}. 
\end{equation}
\end{definition}
In other words, $\ind_\O(x)$ contains all the states in $\O$ that cannot be distinguished from $x$ based on the knowledge of the past output.   To ease the notation, we will omit the mention of $\O$ when $\O = \RR^{n_x}$. As explained in the introduction, in the context of observer design, we are interested in estimating the state modulo its indistinguishable states. For that, we need the preimage of a map $T$ to describe exactly the indistinguishable sets. This leads to the following definition.

\begin{definition}[Characterizing indistinguishability]
\label{def_charac_indist}
A map $T : \mathbb{R}^{n_x} \rightarrow \mathbb{R}^{n_z}$ is said to characterize the backward-distinguishable points in $\O$ for \eqref{sys} if,   for each $(x_a,x_b)\in \O \times \O$, we have
    \begin{equation}
    \label{eq_injectivity_indistinguishable}
       T(x_a) = T(x_b)
       \quad \Longleftrightarrow
       \quad
       x_a\indist x_b.
    \end{equation} 
\end{definition}

The aforementioned notions focus on backward-indistinguishable   \textit{points}.  Next, we define the notion of indistinguishable  \textit{solutions}.

\begin{definition}[Indistinguishable solutions]
Two solutions $\phi_a : I_a\subset \RR\to \RR^{n_x}$ and $\phi_b : I_b\subset \RR\to \RR^{n_x}$ to \eqref{sys} are indistinguishable if
$$ h(\phi_a(t)) = h(\phi_b(t)) \qquad \forall t \in  I_a \cap I_b. $$
\end{definition}

Note that when two points are backward-indistinguishable, then the two maximal backward solutions to \eqref{sys} initialized at those points are indistinguishable. However, backward indistinguishability says nothing about what happens to those solutions in positive time.
Conversely, two indistinguishable solutions verify 
\begin{align}  \label{eqdistus}
\phi_a(t) \indist \phi_b(t) \quad \forall t \in I_a \cap I_b 
\end{align} 
provided that they are \textit{maximal} in backward time.
Note though that, for analytic systems,   the solutions are analytic in time,  and therefore, equality of outputs on an arbitrarily small open subset of  $I_a \cap I_b$, implies equality of outputs on $I_a \cap I_b$.  Hence,  when \eqref{sys} is analytic, two solutions $\phi_a$ and $\phi_b$, with $I_a \cap I_b$ nonempty, are indistinguishable
if and only if there exists $t_o \in  I_a \cap I_b$ such that $\phi_a(t_o) \indist \phi_b(t_o)$. Actually, in this case, \eqref{eqdistus} holds, and the forward and the backward indistinguishability of points are equivalent.

\subsection{Continuity in Set-Valued and $p$-Valued Maps}

In order to study the regularity and convergence of set-valued maps, the space of subsets of $\RR^{n_x}$ is endowed with the Hausdorff distance as defined next.
\begin{definition} [Hausdorff distance]
Given two subsets $\X_a$ and $\X_b$ of $\RR^{n_x}$, the Hausdorff distance is defined as
$$d_\H(\X_a,\X_b) := \max \big\{ \delta(\X_a,\X_b), \delta(\X_b,\X_a)\big\}$$
where 
$$
\delta(\X_a,\X_b) := \sup_{x_a \in \X_a}{d(x_a,\X_b)}  = \sup_{x_a \in \X_a}\inf_{x_b \in \X_b} {d(x_a,x_b)} \ .
$$
\end{definition}

\begin{remark}
Note that, for closed sets, we have
\begin{subequations}
\begin{align}
    \delta(\X_a,\X_b) = 0 & \Longleftrightarrow \X_a\subseteq \X_b \label{eq_delta_inclusion}\\
    d_\H(\X_a,\X_b)=0 & \Longleftrightarrow \X_a=\X_b 
\end{align}
\end{subequations}
\end{remark}

The continuity of a set-valued map $F : K \rightrightarrows \mathbb{R}^{n_x}$ in the sense of the Hausdorff distance contains the following two properties (see \cite[Chapter 1]{aubinDifferentialInclusionsSetValued1984a}).

\begin{definition} [Upper semicontinuity]
A set-valued map $F : K \rightrightarrows \mathbb{R}^{n_x}$ is said to be upper semicontinuous at $z^\star\in K$ if for any open neighbourhood $V$ containing $F(z^\star)$ there exists a neighbourhood $W$ of $z^\star$ such that for all $z \in W$, $F(z) \subseteq V$.
This is equivalent to
\begin{equation}
\label{usc_eq}
\lim_{z\to z^\star}  \delta(F(z),F(z^\star))=0.
\end{equation}
\end{definition}

\begin{definition} [Lower semicontinuity]
A set-valued map $F : K \rightrightarrows \mathbb{R}^{n_x}$ is said to be lower semicontinuous at $z^\star\in K$ if for any $x \in F(z^\star)$ and any neighborhood $V$ of $x$, there exists a neighborhood $W$ of $z^\star$ such that for all $z\in W$, $V \cap F(z) \neq \emptyset$.
This is equivalent to
\begin{equation}
\label{lsc_eq}
\lim_{z\to z^\star} \delta(F(z^\star),F(z))=0.
\end{equation}
\end{definition}


Following the standard topology, the Hausdorff distance can also be used to define stronger continuity property such as Lipschitzness.

\begin{definition} [(Local) Lipschitzness] \label{deflip}
The set-valued map $F : K \rightrightarrows \mathbb{R}^{n_x}$, with $K \subset \mathbb{R}^m$, is said to be \textit{locally Lipschitz} if, for each $x\in K$, there exists a neighborhood $U$ of $x$ and a scalar $k>0$ such that, 
for all $(x_a,x_b) \in (U \cap K) \times (U \cap K)$, 
\begin{align} \label{eq.lipset}
d_\H(F(x_a), F(x_b)) \leq  k |x_a-x_b| . 
 \end{align}
Besides, $F$ is said to be \textit{Lipschitz} if there exists $k>0$ such that \eqref{eq.lipset} holds for all $(x_a,x_b)\in K\times K$.
\end{definition}


Note that a specific and adapted distance may be used on $\A_p(\mathbb{R}^{n_x})$ as defined next.

\begin{definition}[Distance on $\A_p(\mathbb{R}^{n_x})$]
For $S_a \in \A_p(\mathbb{R}^{n_x})$ and $S_b \in \A_p(\mathbb{R}^{n_x})$, we let the distance  
$$
\G(S_a,S_b) := \min_{\sigma \in \P_\nind} \max_{i\in \{1,\ldots,\nind\} } |s_{a,i} - s_{b,\sigma(i)}|
$$
where $\P_\nind$ is the set of permutations of $p$ elements, and $S_a, S_b$ are in $\A_p(\RR^{n_x})$ with elements $\{s_{a,i}\}_{i=1}^p$ and $\{s_{b,i}\}_{i=1}^p$,  respectively. 
\end{definition}

Continuity and (local) Lipschitzness of Almgren $p$-valued maps $\tilde{F} : \Z \subset \mathbb{R}^{n_z} \rightarrow \A_p(\mathbb{R}^{n_x})$, then,  naturally follows with the Hausdorff distance is replaced by the distance $\G$.

\section{Problem Statement} \label{Sec.III}

According to the KKL methodology, we are supposed to find a continuously differentiable map $T: \RR^{n_x} \to \RR^{n_z}$ transforming the dynamics \eqref{sys} into \eqref{z_system};  namely,  solution to \eqref{eq_PDE}. Without observability/distinguishability assumptions,  we cannot hope to prove injectivity of $T$ on $\X$.
However,  \cite[Theorem 1]{PaulineNolcos} showed  the existence of a map $T$ solution to \eqref{eq_PDE} that characterizes the backward-distinguishable states according to Definition \ref{Def-dis-dis}.

\begin{theorem} 
\label{thm_info}
Assume the existence of an open bounded set $\O$,  containing 
$\X$, that is backward invariant for \eqref{sys}.  
Denote $n_o := 2n_x+1$ and $n_z:=n_on_y$.
Then,  there exists $\rho>0$ such that,  for almost any pair $(A_o,B_o) \in \RR^{n_o\times n_o} \times \RR^{n_o}$,  with $A_o+\rho I_{n_o}$ Hurwitz, there exists $T:\RR^{n_x}\to \RR^{n_z}$ continuously differentiable,  characterizing the backward-distinguishable states,  and verifying \eqref{eq_PDE} with  
$
(A = I_{n_y} \otimes A_o,
B = I_{n_y} \otimes B_o)$.

\end{theorem}

\begin{remark}
Note that it is always possible to make the set $\O$ backward invariant by modifying $f$ outside of $\O$.  However,  by doing so,  \eqref{eq_injectivity_indistinguishable} holds for a modified system.  Hence,  $T(x_a)=T(x_b)$ ensures equality of the outputs of the original system only as long as the backward solutions $t\mapsto(\phi_a(t),\phi_b(t))$ from $(x_a,x_b)$ remain in the set where $f$ has not been modified.  
\end{remark}

\begin{remark}
\label{rem_uniqueness}
According to \cite[Proposition 2.1]{BriAndBerSer}, the solution $T$ to \eqref{eq_PDE} is unique on the backward-invariant set $\O$. So, if we find $T$ verifying \eqref{eq_PDE} everywhere on $\O$, for some $n_z \geq n_x$ and for $(A,B)$ picked according to Theorem \ref{thm_info}, then $T$ characterizes the backward-distinguishable states.
\end{remark}

Now, given $t\mapsto x(t)$ solution to \eqref{sys} generating an output $t\mapsto y(t)$,  we know, using \eqref{eq_PDE} that, any solution $t\mapsto z(t)$ to \eqref{z_system}, subject to output $y$, verifies 
\begin{equation}
    \label{eq_conv_KKL_z}
    \lim_{t\to +\infty} |z(t)-T(x(t))|=0.
\end{equation}
Consider the set-valued preimage map $T^-:T(\X) \rightrightarrows \X$ defined by 
\begin{equation} \label{eq_defT-}
T^-(z) := \{ x \in \X : T(x) = z \}.  
\end{equation}
Since $T$ characterizes the backward-indistinguishable states,  in the sense that \eqref{def_charac_indist} holds, we know that 
\begin{align} \label{eqdisdis}
T^-(T(x)) =  \ind_\X(x) \qquad  \forall x \in \X. 
\end{align}
Exploiting \eqref{eq_conv_KKL_z}, it becomes tempting to use $T^-(z(t))$ as an estimate of the indistinguishable set $\ind_\X(x(t))$ at time $t$. However, $z(t)$ is not guaranteed to be in the image set $T(\X)$, where $T^-$ is defined. Therefore,  we  need to find a set-valued extension $\Tinv : \mathbb{R}^{n_z} \rightrightarrows \mathbb{R}^{n_x}$, having the same regularity as $T^-$, and verifying
\begin{equation} \label{eq_TinvT}
    \Tinv(T(x)) =  T^-(T(x))  \qquad \forall x \in \X . 
\end{equation}
From there, one may compute online the set $\Tinv(z(t))$, or at least a continuous selection $\hat{x}(t)\in \Tinv(z(t))$, for instance through optimization schemes. This leads to the following questions:
\begin{enumerate}
    \item Under which conditions does $t \mapsto \Tinv(z(t))$ asymptotically converge to $t \mapsto \ind_\X(x(t))$ in the Hausdorff sense ?
    \item Can we ensure that a continuous selection $t \mapsto \hat{x}(t) \in \Tinv(z(t))$ converges to a solution to \eqref{sys} producing $t \mapsto y(t)$?
\end{enumerate}
Regarding the first question, the answer is positive provided that $\Tinv$ is (Hausdorff) continuous,  which means  upper and lower semicontinuous at the same time. In particular, the latter means that 
$T^-$ must also be (Hausdorff) continuous. 
However, it is shown in \cite{PaulineNolcos} that the map $T^-$ is, in general, only upper semicontinuous. Indeed, the following example shows a case where $T^-$ is not lower semicontinuous and, therefore, it cannot be continuous.  
\begin{example} \label{exp1}
Consider  the smooth map $T : \X :=  [-\pi, \pi]\rightarrow [-1,1] \times [-1,1]$ given by
$$ T(x) := [\sin(2x) ~~ \sin(x)]^\top.   $$
Note that $T^-(0,0) = \{-\pi,0,\pi\}$; whereas, for each $z \in T([-\pi,\pi]) \backslash \{ 0 \}$,  $T^-(z)$ contains only two elements.   Hence, $T^-$ is not lower semicontinuous at $z= 0$.  
\end{example}

As a consequence,  the continuity and convergence of $t \mapsto \Tinv(z(t))$ to  $t \mapsto \ind_\X(x(t))$ is not guaranteed without the extra assumptions we investigate in the next section.  

Before concluding this section, we revisit the example proposed in \cite{PaulineNolcos}, for which, the aforementioned questions have positive answers. 

\begin{example}[\cite{PaulineNolcos}]\label{ex:nolcos}
Consider system \eqref{sys} with
\begin{equation}
\label{num_sys}
\begin{aligned} 
f(x) & := \begin{bmatrix} 
       x_2 + x_1(1-(x_1^2+x_2^2)) \\  -x_1 + x_2(1-(x_1^2+x_2^2)) 
     \end{bmatrix}, 
     \\[0.5em] 
     h(x) & := [x_1^2 - x_2^2 ~~ 2 x_1 x_2 ]^\top.  
 \end{aligned}
  \end{equation}
For any $x\in \RR^2$, $\ind(x) = \{x, -x\}$. 
Next, using the KKL toolbox \cite{KKLtoolbox}, a numerical solution to \eqref{eq_PDE} defined on 
$\X := \B(0,1.7)$ and for $n_z=(n_x+1)n_y=6$ (instead of $n_z=(2n_x+1)n_y=10$ as in Theorem \ref{thm_info}) was obtained,  after making some open set $\mathcal{O}$ backward-invariant. Furthermore, it was numerically checked that $T$ characterizes the backward-distinguishable states in the sense that
$$
T^-(T(x)) = \{x,-x\} \qquad \forall x\in \X.
$$
Moreover, by online solving the optimization algorithm 
$$ \hat{x}(t) := \argmin_{x_s\in \X} |z(t)-T(x_s)|, 
$$
for $z$ solution to \eqref{z_system} subject to an output $y$, 
we found that the obtained function $t \mapsto \hat{x}(t)$ converges asymptotically to the set-valued map $t \mapsto \{x(t),-x(t)\}$ formed by the two solutions to \eqref{sys} generating $y$. However, this does not imply that $t \mapsto \hat{x}(t)$ tends to a solution to \eqref{sys}, unless continuity of $t \mapsto \hat{x}(t)$ is enforced in the optimization problem. In which case, $t \mapsto \hat{x}(t)$ either tends to $t\mapsto x(t)$ or to $t\mapsto -x(t)$; namely, to one of the two indistinguishable solutions. The goal of the paper is to prove these facts observed in simulation.
\end{example}

\section{Lipschitzness of $T^-$ and existence of Lipschitz extension} \label{Sec.IV}

In this section,  we provide assumptions, under which, the set-valued map $T^-$  is Lipschitz (and therefore continuous) and admits a Lipschitz extension $\Tinv$.   More precisely,  we consider the following two assumptions. 

\begin{assumption}[Constant cardinality] \label{ass1a}
There exists $p\in \NN_{>0}$ such that, for each $z\in T(\X)$, $\card(T^{-}(z))=p$.
\end{assumption}

If $T$ characterizes the backward-indistinguishable states as guaranteed in Theorem \ref{thm_info}, namely, \eqref{def_charac_indist} holds, we know that  \eqref{eqdisdis} holds. 
Therefore, Assumption \ref{ass1a} holds if and only if 
$$ \card(\ind_\X(x))=p \qquad 
\forall x\in \X. $$ 
In other words, any state in $\X$ admits exactly $p-1$ other backward-indistinguishable states in $\X$ for \eqref{sys}. For instance, in Example \ref{ex:nolcos}, Assumption \ref{ass1a} holds on any compact subset of $\RR^2\setminus \{0\}$ with $p=2$, but it does not hold on $\X$ since $0\in \X$ and $\card(T^-(T(0))=\card(\ind_\X(0)))=1$.

\begin{assumption}[Full rank] \label{ass1b} 
 For each $x \in \X$, the jacobian matrix $ \frac{\partial T}{\partial x}(x)$ is full rank.   
 \end{assumption}   
 
\begin{figure}
\label{fig:jacobian}
\includegraphics[width=8cm]{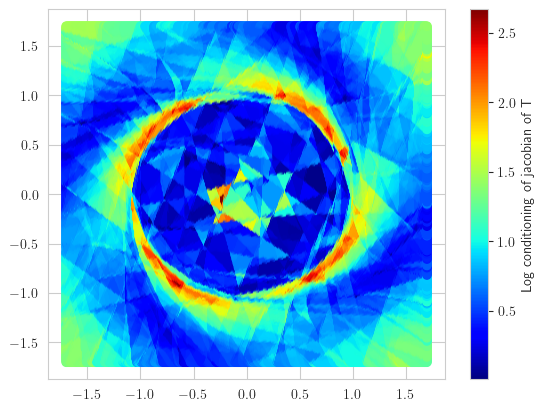}
\caption{Condition number of the jacobian of $T$ in Example \ref{ex:nolcos}}
\end{figure}

 The ``physical'' interpretation of Assumption \ref{ass1b} is less straight-forward. Indeed, when the map $T$ solving \eqref{eq_PDE} is known explicitly or learned numerically as in \cite{KKLtoolbox}, the conditioning of the jacobian matrix of $T$ can be evaluated on a numerical grid. For instance, the jacobian of $T$ obtained in \cite{PaulineNolcos}, for the system in Example \ref{ex:nolcos}, can be checked to be full-rank on $\X$ with a conditioning smaller than $10^3$; see Figure \ref{fig:jacobian}. Furthermore, in Section \ref{sec_diff_obs}, we  exhibit a link between the rank of $T$ and the rank of the differential observability map in \eqref{eq_def_H}, provided that the eigenvalues of $A$ have a sufficiently large real part.
 The following lemma discusses a link between Assumptions \ref{ass1a} and \ref{ass1b}.
\begin{lemma}
Let $T: \O \subset\RR^{n_x} \to \RR^{n_z}$ be a continuously differentiable map, and let $\X \subset \O$ be compact such that Assumption \ref{ass1b} holds. Then, $\card (T^-)$, with $T^-$ defined as in \eqref{eq_defT-}, is finite.
\end{lemma}
\begin{pf}
Let $x_o\in \X$ and $z_o=T(x_o)$. Under Assumption \ref{ass1b}, for any $x_o'\in T^-(z_o)\subseteq \X$, there exists a neighborhood $\BB(x_o',\varepsilon)$ such that $T:\BB(x_o',\varepsilon)\to \RR^{n_z}$ is injective. Therefore, $x_o'$ is the unique preimage of $z_o$ in $\BB(x_o',\varepsilon)$, i.e., the preimages of $z_o$ are isolated. Since $\X$ is compact, $T^-(z_o)$ is finite.  
\end{pf}

We conclude from this lemma that the contribution of Assumption \ref{ass1a} over \ref{ass1b} is not in the fact that $\card(T^-)$ is finite, but in the fact that it is constant. For instance, in Example \ref{exp1}, the map $T$ verifies Assumption \ref{ass1b}, and yet does not show constant cardinality of $T^-$. 
This is also the case of Example \ref{ex:nolcos}, where Assumption \ref{ass1b} holds and yet the cardinality drops at zero.  In fact, this is because $T$ is not a local \textit{homeomorphism} on $\X$.
On the other hand, a map $T$ verifying Assumption \ref{ass1a} does not necessarily verify Assumption \ref{ass1b} (take, for instance, the map $T : \mathbb{R} \rightarrow \mathbb{R}$ with $T(x) := x^3$). 

Under Assumptions \ref{ass1a} and \ref{ass1b}, we establish  Lipschitzness of $T^-$ on $T(\X)$. Furthermore, thanks to Almgren theory \cite{Alm}, we establish the existence of a set-valued map $\Tinv : \mathbb{R}^{n_z} \rightrightarrows \mathbb{R}^{n_x}$ that is a 
Lipschitz extension of $T^-$ to $\mathbb{R}^{n_z}$.  

\begin{theorem} \label{theo_Lips_Tinv}
Consider a $C^1$ map $T : \O \subset \RR^{n_x} \to \RR^{n_z}$, with $\O$ an open set containing $\X$ and $n_z \geq n_x$, such that Assumptions \ref{ass1a}  and \ref{ass1b} hold. 
Then, there exists a Lipschitz set-valued map $\Tinv : \RR^{n_z} \rightrightarrows \RR^{n_x}$ such that
\begin{align} \label{eqTinvrestr}
\Tinv(z) = T^-(z) \qquad \forall z\in T(\X)
\end{align}
and for any $z\in \RR^{n_z}$, $\Tinv(z)$ contains at most 
$p$ elements.  
\end{theorem}

\begin{pf}
Consider $z_o \in T(\X)$, and according to Assumption \ref{ass1a},  there exist $(x_1,x_2, \ldots, x_p) \in \X^p$ such that 
$$
T^-_\X(z_o) = \{x_1,\ldots,x_p\} .
$$
Applying the  constant rank theorem in Lemma \ref{lem1},  we conclude that,  for each $i\in \{1,2,\ldots, p\}$,   there exist $\varepsilon_i,  ~ \sigma_i > 0$ and two diffeomorphisms 
$\phi_i: \B(x_i, \varepsilon_i)  \to \RR^{n_x}$ and 
$\psi_i : \B(z_o, \sigma_i)  \to \RR^{n_z}$  such that 
$$
 \psi_i \circ T(x) = (\phi_i(x), 0, \ldots , 0) \qquad \forall x \in \B(x_i, \varepsilon_i).
$$
Now, denote by $\text{Proj} :\RR^{n_z} \to \RR^{n_x}$ the projection on the the first $n_x$ coordinates, we let the $C^1$ map $g_i := \phi_i^{-1} \circ \text{Proj} \circ \psi_i: \B(z_o, \sigma_i) \to \RR^{n_x}$  verifying 
$$ 
g_i \circ T(x) = x \qquad \forall x \in  \B(x_i, \varepsilon_i). 
$$ 
In other words,  $g_i$ is a $C^1$ local left-inverse of $T$ around $x_i$.   We now prove that there exists $\sigma_o > 0$ such that
\begin{equation}
\label{eq_union_gi}
    T^-(z) = \bigcup_{i=1}^\nind g_i(z) \qquad \forall z\in \B(z_o,\sigma_o) \cap T(\X).
\end{equation}
To do so,  we pick a sequence $(z_k)_{k\in \NN}$ in $T(\X)$ converging to $z_o$.   
Since 
$$ T^-(z_o) := \{x_1, x_2,\ldots,x_p\} \subset \cup_{i=1}^p \B(x_i, \varepsilon_i),  $$ 
and since $T^-$ is upper semicontinuous acording 
to \cite{PaulineNolcos},  it follows that there exists $k^*>0$ such that,  for each $k\geq k^*$,  
$$ T^-(z_k)\subset \cup_{i=1}^p \B(x_i, \varepsilon_i).   $$ 
Besides, from Assumption \ref{ass1a},  we have
$$ T^-(z_k)=\{x_{k,1},  x_{k,2},\ldots,x_{k,p}\} \qquad \forall k\in \NN,  $$ 
with $ x_{k,l} \neq x_{k,j} $ whenever $ l\neq j$.
As a result,  for each $k\geq k^*$,  and for each $j\in \{1,2,\ldots,p\}$,  there exists $i_{kj}\in \{1,2,\ldots,p\}$ such that
$$
z_k = T(x_{k,j}),  \quad \text{and} \quad x_{k,j} \in \B(x_{i_{kj}} ,  \varepsilon_{i_{kj}}).
$$
Thus,  $x_{k,j} = g_{i_{kj}}(z_k)$,  so that
$$T^-(z_k)\subseteq \cup_{i=1}^p g_{i}(z_k) \qquad \forall k\geq k^*.  $$
The latter establishes the existence of a neighborhood $W_o$ of $z_o$ such that 
$$ T^-(z) \subseteq \bigcup_{i=1}^\nind g_i(z) \qquad  \forall z \in W_o \cap T(\X).  $$ 
However, by Assumption \ref{ass1a},  the cardinal of $T^-(z)$ is $p$ for $z \in T(\X)\cap W_o$,  so that necessarily it is an equality and \eqref{eq_union_gi} follows. 
Then, we note that the set-valued map $g : W_o \rightrightarrows \RR^{n_x}$ defined by  
$g(z) = \bigcup_{i=1}^\nind g_i(z)$
can be identified with an Almgren $\nind$-valued map 
$\tilde{g}$ on $W_o \cap T(\X)$  
Thus, to $g$ and $T^-$,  we associate  the Almgren
$\nind$-valued maps $\tilde{g}$ and $\tilde{T}^-$ defined on $W_o$ and $W_o \cap T(\X)$,   respectively.  
  
Since each $g_i$,  $i \in \{1,2,...,p\}$,  is $C^1$,  it is  Lipschitz on $W_o$ with a Lipschitz constant $L_i > 0$,  and after \eqref{eq_union_gi},  we conclude that,  for each $z_a,z_b\in W_o  \cap T(\X)$,  we have 
\begin{align*}
    \G(\tilde{T}^-&(z_a),\tilde{T}^-(z_b))= \min_{\sigma \in \P_\nind} \max_{i\in \{1,\ldots,\nind\}} |g_{i}(z_a) - g_{\sigma(i)}(z_b)|   \\
 &  \leq \max_{i\in \{1,\ldots,\nind\}} |g_{i}(z_a) - g_{i}(z_b)| \leq \max_{i\in \{1,\ldots,\nind\}} L_i \ |z_a-z_b|. 
\end{align*}
This being true on a neighborhood of $z_o$,  for any $z_o \in T(\X)$,  we conclude that  $\tilde{T}^-$ is locally Lipschitz on $T(\X)$.    Now,  since $T(\X)$ is compact under the continuity of $T$,  we show that $\tilde{T}^-$ is (globally) Lipschitz on $T(\X)$.    Indeed,   assume $\tilde{T}^-$ is not Lipschitz on  $T(\X)$. Then, for every $k\in \NN$, there exists $z_{a,k},z_{b,k}$ in $T(\X)$ such that $ \G(\tilde{T}^-(z_{a,k}),\tilde{T}^-(z_{b,k}))\geq k |z_{a,k}-z_{b,k}|$.  By compactness, we can assume that the sequence $\{z_{a,k},z_{b,k} \}^{\infty}_{k=1}$ is convergent to $(z_a^*,z_b^*)\in T(\X) \times T(\X)$.   Since 
$\tilde{T}^-(z_{a,k})$ and $\tilde{T}^-(z_{a,k})$ are bounded, necessarily  $|z_{a,k}-z_{b,k}|$ converges to $0$, and, thus, $z_{a}^*=z_b^*:=z^*$. Therefore, for sufficiently large $k$, $z_{a,k},z_{b,k}$ must belong to a neighborhood of $z^*$ where $\tilde{T}^-$ is Lipschitz.  This yields to a contradiction when  $k$ is sufficiently large. 

Next,  by using Almgren extension theorem; see Lemma \ref{lem2} in the Appendix,
we conclude the existence of a Lipschitz Almgren $p$-valued map $\tilde{T}^{\rm inv} : \RR^{n_z} \to \A_p(\RR^{n_x})$ (with multiplicity allowed) such that $\tilde{T}^{\rm inv}$ agrees with $\tilde{T}^-$ on $T(\X)$.  
In particular, there exists $L > 0$ such that
$$
\G(\tilde{T}^{\rm inv}(z_a),\tilde{T}^{\rm inv}(z_b))\leq L |z_a-z_b| \quad \forall (z_a,z_b)\in \RR^{n_z} \times \RR^{n_z} \ .
$$
Then, forgetting about multiplicity and transforming the $p$-valued map into a set-valued one,  we build $\Tinv : \RR^{n_z} \rightrightarrows \RR^{n_x}$ such that 
$$
d_\H(\Tinv(z_a),\Tinv(z_b))\leq L |z_a-z_b| \quad \forall (z_a,z_b)\in \RR^{n_z} \times \RR^{n_z}.
$$
Indeed,
\begin{multline*}
    d_\H(\Tinv(z_a),\Tinv(z_b))\\= \max \{ \max_{i} \min_{j} |\tilde{T}_i^{\rm inv}(z_a) - \tilde{T}_j^{\rm inv}(z_b) | ,\\ \max_{j} \min_{i} |\tilde{T}_i^{\rm inv}(z_a) - \tilde{T}_j^{\rm inv}(z_b) | \}
\end{multline*}
and for any permutation $\sigma$ in $\P_p$, 
\begin{align*}
   \max_{i} \min_{j} |\tilde{T}_i^{\rm inv}(z_a) - \tilde{T}_j^{\rm inv}(z_b) |
   &  \\ & \hspace{-0.8cm} \leq \max_{i}  |\tilde{T}_i^{\rm inv}(z_a) - \tilde{T}_{\sigma(i)}^{\rm inv}(z_b) | \\
   \max_{j} \min_{i} |\tilde{T}_i^{\rm inv}(z_a) - \tilde{T}_j^{\rm inv}(z_b) |& \\ & \hspace{-0.8cm} \leq \max_{j}  |\tilde{T}_{\sigma^{-1}(j)}^{\rm inv}(z_a) - \tilde{T}_{j}^{\rm inv}(z_b) |, 
\end{align*}
which gives the result by an appropriate re-indexing.
\end{pf}

\begin{example}
Consider the map $T$ introduced in Example \ref{exp1}. We already showed that $T^-$ is only upper semicontinuous.   We can see that $\frac{\partial T}{\partial x}(x)$ is full rank for all $x \in \X := [-\pi,\pi]$.  Hence,  Assumption \ref{ass1b} is verified.  However,  Assumption \ref{ass1a} is not verified on $T(\X)$ but verified on any compact set contained in $T(\X) \backslash \{ 0 \}$.   Thus,  $T^-$ is Lipschitz on any compact subset contained in $T(\X) \backslash \{ 0 \}$.
\end{example}

\begin{remark}
Note that the construction of the extension of $\tilde{T}^{-}$ according to the proof in \cite[Theorem 1.7]{De_Lellis_2011}  is  constructive;  hence,  an algorithm can be deduced to compute,  or at least to approximate,  $\Tinv$.  
\end{remark}

\begin{remark}
Note that the conclusions of Theorem \ref{theo_Lips_Tinv} remain valid if, instead of Assumption \ref{ass1a}, we assume that the cardinality of $T^-$ is constant only on the connected components of $T(\X)$. The cardinality of $\Tinv$ will then not exceed the maximal
cardinality of $T^-$. Such a relaxation will not help our study since,  using \eqref{eq_conv_KKL_z}, a solution $z$ to \eqref{z_system}, subject to the system's output $y$, converges to $T(\X)$, and, since $z$ is continuous, then it necessarily converges to a connected component in $T(\X)$, on which, the cardinality of $T^-$ is constant under both assumptions. 
\end{remark}

\section{Convergence of the KKL Observer} \label{sec_conv_KKL}

\subsection{Set-valued convergence to indistinguishable set}

We start this section by establishing the following direct consequence of Theorem \ref{theo_Lips_Tinv}, which holds under  the following assumption.
\begin{assumption} \label{assdist}
The map $T: \X  \to \RR^{n_z}$  characterizes the backward-distinguishable points in $\mathcal{X}$ for \eqref{sys} according to Definition \ref{def_charac_indist}.
\end{assumption}
The existence of a map $T$ verifying Assumption \ref{assdist} is guaranteed by Theorem \ref{thm_info}. Then, as noticed above, Assumption \ref{ass1a}  holds if and only if any point in $\X$ admits exactly $p-1$ other backward indistinguishable states in $\X$.

\begin{theorem} \label{theo_hausdorff_conv}
Consider a $C^1$ map $T : \O \subset \RR^{n_x} \to \RR^{n_z}$, with $\O$ open and $n_z \geq n_x$, verifying \eqref{eq_PDE}, for a given pair $(A,B)$ and $\X$ contained in $\O$, and satisfying Assumptions \ref{ass1a}, \ref{ass1b}, and \ref{assdist}. Let $\Tinv:\RR^{n_z} \to \RR^{n_x}$ be the Lipschitz extension given in Theorem \ref{theo_Lips_Tinv}.  Then, there exist $\rho,\lambda>0$ such that, for any solution 
$x:\RR_{\geq 0} \to \RR^{n_x}$ to \eqref{sys}  remaining in $\X$, and for any solution $z:\RR_{\geq 0} \to \RR^{n_z}$ to \eqref{z_system} subject to $y=h(x)$,  we have
\begin{multline}
\label{eq_hausdorff_exp_conv}
d_\H(\Tinv(z(t)),\ind_\X(x(t))) \\ \leq \rho e^{-\lambda t} |z(0)-T(x(0))|  \quad \forall t\geq 0.  
\end{multline}
\end{theorem}
Therefore, if the system solution $x$ eventually remains in $\X$, then the set-valued estimate $\Tinv(z(t))$ asymptotically converges in the Hausdorff sense to the indistinguishable set $\ind_\X(x(t))$.

\begin{pf}
According to Theorem \ref{theo_Lips_Tinv}, $\Tinv$ is (globally) Lipschitz and verifies 
$$ \Tinv(T(x))  = T^- (T(x))  \qquad 
\forall x \in \X. $$
Furthermore, using Assumption \ref{assdist}, we conclude that 
$$ T^- (T(x)) = \ind_\X(x) \qquad 
\forall x \in \X. $$
Therefore, there exists $L>0$ such that
$$
d_\H(\Tinv(z(t)),\ind_\X(x(t)))\leq L |z(t)-T(x(t))| 
$$
for all $t\geq 0$ provided that $x(t)\in \X$. Hence, the proof is completed  using the fact that the estimation error $\tilde{z}(t):=z(t) - T(x(t))$ in the $z$-coordinates verifies the exponentially stable dynamics $\dot{\tilde{z}}=A\tilde{z}$.
\end{pf}

\begin{remark}[ISS set-valued observer]
If the available output $y$ fed into \eqref{z_system} is noisy, i.e., $y=h(x)+\nu$, then, the error $\tilde{z}(t) := z(t)-T(x(t))$ in the $z$-coordinates evolves according to $\dot{\tilde{z}}=A\tilde{z} + B\nu$ instead of $\dot{\tilde{z}}=A\tilde{z}$ and there exist $\rho_1,\rho_2,\lambda>0$ (depending on the pair $(A,B)$ only) such that
\begin{multline} \label{eq_ISS}
    d_\H(\Tinv(z(t)),\ind_\X(x(t)))\\\leq L \left(\rho_1 e^{-\lambda t} |z(0)-T(x(0))| + \rho_2 \sup_{s\in [0,t]} |\nu(s)| \right)
\end{multline}
for all $t\geq 0$ provided that $x(t) \in \X$. Besides, if $x(0)\in \X$ and $z(0)$ is chosen equal to $T(\hat{x}_0)$ for some initial guess $\hat{x}_0\in \X$, then $|z(0)-T(x(0))|\leq L_T|\hat{x}_0-x(0)|$, where $L_T$ is the Lipschitz constant of the $C^1$ map $T$ on the compact set $\X$. It follows from \eqref{eq_ISS} that the KKL observer exhibits a set-valued ISS property with respect to measurement noise, in the Hausdorff sense.
\end{remark}

\begin{remark}
In the case where $t\mapsto x(t)$ does not remain in $\X$, the map $t \mapsto \Tinv(z(t))$ is  guaranteed to approach $t \mapsto \ind_\X(x(t))$ only during the time intervals on which the solution $x$ is within the set $\X$. 
\end{remark}


\subsection{Convergence of a continuous selection to a solution}

In practice,  one can compute a \textit{continuous} selection $t \mapsto \hat{x}(t) \in \Tinv(z(t))$,  for instance, through an optimization algorithm. Hence, it is interesting to know whether $\hat{x}$ converges to a solution generating $y$ or not. To answer this question, we make the following assumption on the solutions generating a same output $y$.


\begin{assumption} \label{asssy}
Given $y:\RR_{\geq 0}\to \RR^{n_y}$, there exist $\bar{t}>0$ and at least $p \in \mathbb{N}_{>0}$ solutions $ \{ x_i \}^p_{i=1}: [\bar{t},+\infty) \to \X$ to \eqref{sys} such that, for all $i \in 
\{1,\ldots,p\}$, we have 
\begin{itemize}
    \item $h(x_i(t))=y(t)$ for all $t\geq \bar{t}$,
    \item for all $j\in \{1,\ldots,p\}$ with $j\neq i$, $x_i-x_j$ does not asymptotically converge to zero. 
\end{itemize}
\end{assumption}

Given a solution $z$ to \eqref{z_system} subject to $y$, we would like to relate the continuous selection $t \mapsto \hat{x}(t) \in \Tinv(z(t))$  
and the solutions $x_i$ to \eqref{sys} generating $y$. For that, we rely on the following \textit{Splitting Lemma}
first proved in \cite{Banach1934}, see also  \cite[Theorem 3.1.]{nvalued2018} for more details.

\begin{lemma}[Splitting Lemma] \label{lem3}
Consider a continuous set-valued map 
$T^- : \C \subseteq \mathbb{R}^{n_z} \rightrightarrows \mathbb{R}^{n_x} $ such that $\C$ is simply connected and there exists $p\in \NN_{>0}$ such that, for each $z\in \C$, $\card(T^-(z))=p$.
Then, there exist $p$ continuous (single-valued) functions $\{ g_i \}^p_{i=1}  : \C \rightarrow \mathbb{R}^{n_x}$  such that 
\begin{equation}\label{eq_TminSplit}
    T^-(z) = \bigcup^p_{i=1}  \{g_i(z)\} \qquad \forall z \in \C.
\end{equation}
The map $T^-$ is then said to be \textit{split} on $\C$.
\end{lemma}

\begin{remark}
Without any assumption on its domain $\C$ (such as simple connectedness), a continuous $p$-valued map is not necessarily split in the general case where $n_z$ and $n_x$ are both greater than one \cite{staecker2021partitions}; see the following example. 
\end{remark}

\begin{example}
Let $\C \subset \mathbb{R}^2$ be the unite circle; namely, 
$$ \C := \{ z \in \mathbb{R}^2 : |z| = 1 \}. $$
We parameterize $\C$ by the variable $t \in [0,1)$, and we let the map 
$\phi : [0,1) \rightarrow \C$ defined as
$$ \phi(t) := ( \sin(2\pi t), \cos(2\pi t) ). $$
Note that $\phi$ is a diffeomorphism and 
$$ \phi^-(z) = \frac{\arg(z)}{2\pi} \qquad \forall z \in \C,  $$
where $\arg : \C \rightarrow [0,2 \pi)$ is the function associating each $z \in \C$ to its argument.

Consider the set-valued map 
$H : [0,1) \rightrightarrows \mathbb{R}^2$ given by 
$$ H(t) := \{ f(t), g(t) \},  $$
where $f$, $g : [0,1) \rightarrow \mathbb{R}^2$ are continuous functions
defined as
\begin{align*}
f(t) = (f_1(t),f_2(t)) & := (2-t, 4-t),
\\
g(t) = (g_1(t),g_2(t)) & := (1 + t, 3 + t e^{t-1}).
\end{align*}
We start noting that  
$$ f(t)  \neq g(t) \qquad \forall t \in [0,1). $$
Indeed, $f_1(t) = g_1(t)$ only when $t = 1/2$. However, 
$$ 7 /2 = f_2(1/2) \neq g_2(1/2) = 3 + 1/(2 e^{1/2}). $$  
Hence, $H$ is a continuous set-valued map such that $\card(H(t))=2$ for all $t\in [0,1)$.  
  
At this point, we define the set-valued map $T^- : \C \rightrightarrows \mathbb{R}^2$ given by   
$$ T^-(z) = H(\phi^-(z))  \qquad \forall x \in \C. $$
Clearly, $\card(T^-(z))=2$ for all $z\in \C$. Furthermore, 
since $\phi$ is a diffeomorphism and $\phi^-$ is continuous on $\C \backslash (1,0)$, we conclude that $T^-$ is continuous on $\C \backslash (1,0)$. Actually, observing that
\begin{align*}
\lim_{t \rightarrow 1} f(t) = g(0), \quad \text{and} \quad  \lim_{t \rightarrow 1} g(t) = f(0),
\end{align*}
it follows that $T^-$ is also continuous at $z = (1,0)$.
The latter shows that $T^-$ is $2$-valued and continuous on $\C$. But we next show that $T^-$ is not split using contradiction. Indeed, 
assume that $T^-$ is split, then we can find continuous functions $t^-_1$, $t^-_2 : \C \rightarrow \mathbb{R}^2$ 
such that 
$$ T^-(z) = \{ t^-_1(z), t^-_2(z) \} \qquad z \in \C. $$
At the same time, by definition, we have that 
$$ T^-(z) = \{ f(\phi^-(z)), g(\phi^-(z)) \} \qquad z \in \C. $$
Now, since $f \circ \phi^- $ and $g \circ \phi^- $ are continuous on $\C \backslash (1,0)$ and never intersect, it follows that either 
$$ f(\phi^-(z)) = t^-_1(z) \quad \text{and} \quad g(\phi^-(z)) = t^-_2(z)$$
for all $z \in \C \backslash (1,0)$ or 
$$ f(\phi^-(z)) = t^-_2(z) \quad \text{and} \quad g(\phi^-(z)) = t^-_1(z)$$
for all $z \in \C \backslash (1,0)$. Both cases yield a contradiction since $t^-_1$ and $t^-_2$ are continuous whereas the functions $f \circ \phi^- $ and $g \circ \phi^- $ are discontinuous at $(1,0)$. Therefore, $T^-$ is not split.
\end{example}

In order to apply Lemma \ref{lem3}, given an output  $y$ to \eqref{sys}, we introduce the following assumption.

\begin{assumption} \label{asszz}
There exists a solution $z_*$ to \eqref{z_system} subject to $y$
and a simply connected subset 
 $ \mathcal{C} \subseteq T(\mathcal{X})$ such that $\lim_{t \rightarrow +\infty} |z_*(t)|_{\C} = 0$. 
\end{assumption}

\begin{remark}
Given any output $y$ generated by a solution $x$ to $\eqref{sys}$, any solution $z$ to \eqref{z_system} subject to $y$ converges to $T(\X)$. Therefore, Assumption \ref{asszz}  holds if $T(\mathcal{X})$ is simply connected or if $t\mapsto T(x(t))$ forms (at least after a certain time) an arc that can be continuously reduced to a point in $T(\mathcal{X})$. This holds in particular if $x$ asymptotically converges to a point.
\end{remark}

Assumption \ref{asszz} allows us to guarantee that $T^-$ is split on $\C$. We then show that the same property holds for $\Tinv$ on a neighborhood of $\C$. This allows to show that any continuous selection $t \mapsto \hat{x}(t) \in \Tinv(z(t))$, with $z$ solution to \eqref{z_system} subject to $y$, converges to one of the solutions generating $y$ given by Assumption \ref{asssy}, as stated in Theorem \ref{thm4} below.

But, actually, we can prove that a continuous selection within $t \mapsto \Tinv(z(t))$ converges to a solution  even if $\Tinv$ is not split on a set containing the solutions after a certain time. So we propose the following alternative assumption when Assumption 
\ref{asszz} does not hold. 

\begin{assumption} \label{assz1}
There exists a solution $z_*$ to \eqref{z_system} subject to $y$, for which, there exists $\tau_o > 0$ and $\delta>0$ such that, for each $\tau\geq \tau_o$, $z^*([\tau, \tau+\delta])$ is simply connected.
\end{assumption}

We can then state the following result.

\begin{theorem} \label{thm4}
Consider a $C^1$ map $T : \O \subset \RR^{n_x} \to \RR^{n_z}$, with $\O$ open and $n_z \geq n_x$,
verifying \eqref{eq_PDE}, for a given pair $(A,B)$ and $\X$ contained in $\O$, and Assumptions \ref{ass1a} and \ref{ass1b}. 
Consider an output $y:\RR_{\geq 0}\to \RR^{n_y}$ to \eqref{sys} such that Assumption \ref{asssy} and either Assumption \ref{asszz} or \ref{assz1} hold. 
Then, for any continuous selection $t \mapsto \hat{x}(t)\in \Tinv(z(t))$, where  $z:\RR_{\geq 0}\to \RR^{n_z}$ is a solution to \eqref{z_system} subject to $y$, there exists a solution $\tilde{x}$ to \eqref{sys} generating $y$ such that
\begin{equation}
\label{eq_conv_indist}
\lim_{t \to +\infty} |\hat{x}(t) - \tilde{x}(t)|  = 0.
\end{equation}
\end{theorem}

\begin{pf}
Consider an output $y:\RR_{\geq 0}\to \RR^{n_y}$ to \eqref{sys} and let $z:\RR_{\geq 0}\to \RR^{n_z}$ be a solution to \eqref{z_system} subject to the output $y$.  We start by showing the result under 
Assumption \ref{asszz}.
For that, let $\mathcal{U}_c \subset \mathbb{R}^{n_z}$ be a simply connected compact subset including $\C$ in its interior. Next, we let $\Tinv_u$ be the restriction of $\Tinv$ to $\mathcal{U}_c$. Using Theorem \ref{theo_Lips_Tinv}, we know that $\Tinv_u$ is continuous and its images contain at most $p$ elements. Since $\Tinv$ equals $T^-$ on the compact set $\mathcal{U}_c\cap T(\X)$, whose images contain exactly $p$ elements from Assumption \ref{ass1a}, it follows that, by choosing the boundary of $\mathcal{U}_c$ sufficiently close to $T(\X)$, we can say without loss of generality
\begin{align} \label{eqcardinflat}
\card(\Tinv_u(z)) = p \qquad \forall z \in \mathcal{U}_c. \end{align}
Now, since $\Tinv_u$ is continuous on $\mathcal{U}_c$ simply connected and \eqref{eqcardinflat} holds, 
we conclude, using Lemma \ref{lem3}, the existence of a sequence of continuous functions $\{ g_i \}^p_{i=1}  :  \mathcal{U}_c \rightarrow \mathcal{X}$  such that 
$$ \Tinv_u(z) = \bigcup^p_{i=1}  g_i(z) \qquad \forall z \in \mathcal{U}_c. $$
According to \eqref{eqcardinflat}, we have for all 
$z \in \mathcal{U}_c$ and for all $i\neq i'$, $g_i(z)\neq g_{i'}(z)$. By compactness of $\mathcal{U}_c$, continuity of the maps $g_i$ and finiteness of pairs $i,i'$, there exists $\varepsilon>0$ such that
\begin{equation}
\label{eq_gi_distinct}
    |g_i(z)-g_{i'}(z)|\geq \varepsilon \quad \forall i\neq i', \quad \forall z\in \mathcal{U}_c.
\end{equation}

On the other hand, knowing that the solution $z_*$ in 
Assumption \ref{asszz} converges to $\C$, we know that $z_*(t) \in \mathcal{U}_c$ at least after a certain time.
And knowing that any solution $z$ to \eqref{z_system} subject to $y$  converges to $z_*$, we have by continuity that 
\begin{align} \label{eqlimstart}
\lim_{t \rightarrow \infty} d \left( \hat{x}(t),  \Tinv_u(z_*(t)) \right) = 0. 
\end{align}
Hence, since $\hat{x}$ is continuous,   \eqref{eqlimstart} and \eqref{eq_gi_distinct}  imply  the existence of  $i^* \in \{1,2,...,p\} $ such that 
\begin{align} \label{eqlim2}
\lim_{t \rightarrow \infty} |\hat{x}(t) - g_{i^*}(z_*(t))| = 0.
\end{align}

Now,  we show that each $t \mapsto g_{i}(z_*(t))$, $i \in \{1,2,...,p\}$, must converge to one of the solutions $x_i$ generating $y$ given by Assumption \ref{asssy}. Indeed, let  $\{z_i\}_{i=1}^p:\RR_{\geq 0}\to \RR^{n_z}$ such that 
\begin{equation}
\label{eq_zi}
    z_i(t) = T(x_i(t)) \qquad \forall i \in \{1,2,...,p \}, \quad  \forall t \geq 0.
\end{equation}
According to \eqref{eq_PDE} and the first item of Assumption \ref{asssy}, all the $z_i$s and $z_*$ satisfy \eqref{z_system} subject to $y$. Hence, they converge asymptotically to each other.  Furthermore, according to \eqref{eq_zi} and since $x_i(t)\in \X$ for all $t \geq \bar{t}$, we have
$$ x_i(t) \in  T^-(z_i(t)) \qquad \forall t \geq \bar{t}, \quad  \forall i \in \{1,2,...,p\}. $$ 
As above, using the continuity of $\Tinv$ given by Theorem \ref{theo_Lips_Tinv}, we deduce that
$$
\lim_{t \rightarrow \infty} d \left( x_i(t),  \Tinv_u(z_*(t)) \right) = 0 \qquad \forall i \in \{1,\ldots p\}.
$$
Now, according to the second item of Assumption \ref{asssy} and the fact that $t \mapsto \Tinv_u(z_*(t))$ describes $p$ continuous functions of time not converging to each other according to \eqref{eq_gi_distinct}, we conclude that, for each $i \in \{1,\ldots,p \}$,  there exists a unique $k_i \in \{1,2,...,p\}$ such that $x_{i}$ converges to $ g_{k_i}(z_*)$ asymptotically. Similarly, 
for each $i \in \{1,\ldots,p \}$,
there exists a unique $l_i \in \{1,2,...,p\}$ such that $g_{i}(z_*)$  converges to $x_{l_i}$ asymptotically. From \eqref{eqlim2}, we deduce the result.  


We now suppose Assumption \ref{assz1} holds instead of 
Assumption \ref{asszz}. Let $\mathcal{U} \subset \mathbb{R}^{n_z}$ be a compact subset including $T(\X)$ in its interior and let $\Tinv_u$ be the restriction of $\Tinv$ to $\mathcal{U}$. By Theorem \ref{theo_Lips_Tinv}, $\Tinv_u$ is continuous and its images contain at most $p$ elements. Besides, $\Tinv_u=T^-$ on the compact set $T(\X)$ where its images have exactly $p$ distinct elements according to Assumption \ref{ass1a}.  It follows that, without loss of generality, when $\partial \mathcal{U}$ is sufficiently close to $T(\X)$, then we may assume that
 \begin{align} \label{eqcardinfl2} 
 \card(\Tinv_u(z)) = p \qquad \forall z \in \mathcal{U}. 
 \end{align}
\cite[Proposition 4.1.]{nvalued2018} then shows the existence of $\bar{\epsilon} > 0$ such that
\begin{equation}
\label{eq_defepsilonTinvu}
    \min \{ |x_a-x_b| : x_a \neq x_b \in \Tinv_u(z), ~ z \in \mathcal{U} \}  \geq \bar{\epsilon}. 
\end{equation}

Furthermore, since $y$ is generated by a solution to $\eqref{sys}$ that is eventually in $\X$ according to Assumption \ref{asssy}, any solution to \eqref{z_system} fed with $y$ converges to $T(\X)$ so that  $z_*(t) \in \mathcal{U}$  at least after some time, let's say $\tau_o$ given by Assumption \ref{assz1} without loss of generality.

Finally, all solutions to \eqref{z_system} converge to each other, so in particular $z$ converges to $z_*$ and we have
\begin{align} \label{eqaddednew}
\lim_{t \rightarrow \infty} d \left( \hat{x}(t),  \Tinv_u(z_*(t)) \right) = 0. 
\end{align}
  
Now, under Assumption \ref{assz1}, for each $\tau \geq  \tau_o$,  the set $z_*([\tau , \tau+\delta]) \subset \mathcal{U}$ is simply connected. Hence, using Lemma \ref{lem3}, we conclude the existence of a sequence $\{ g^{\tau}_i \} : z_*([\tau, \tau+\delta]) \rightarrow \X$ such that 
$$ \Tinv_u(z) = \cup^{p}_{i=1} g^{\tau}_i(z) \qquad \forall z \in z_*([\tau, \tau+\delta]). $$
Then, \eqref{eq_defepsilonTinvu} allows us to conclude that 
\begin{align} \label{eqsepar}
|g^{\tau}_i(z) - g^{\tau}_{j\neq i} (z) |  \geq \bar{\epsilon} \quad \forall z \in z_*([\tau,\tau+\delta]), ~~ \forall \tau > \tau_o.   
\end{align}

Let  $\{z_i\}_{i=1}^p:\RR_{\geq 0}\to \RR^{n_z}$ such that 
\begin{equation}
\label{eq_zibis}
    z_i(t) = T(x_i(t)) \qquad \forall i \in \{1,2,...,p \}, \quad  \forall t \geq 0.
\end{equation}
According to \eqref{eq_PDE} and the first item of Assumption \ref{asssy}, all the $z_i$s and $z_*$ satisfy \eqref{z_system} subject to $y$. Hence, they converge asymptotically to each other and are in $\mathcal{U}$ after a certain time. $\Tinv_u$ being continuous on $\mathcal{U}$, and reasoning as above, we conclude that, for each $\epsilon_o \in (0,\bar{\epsilon}/8]$, there exists $\tau_1 \geq \tau_o$ such that, for each $\tau \geq \tau_1$ the following two properties hold:
\begin{itemize}
\item For each $i \in \{1,2,...,p\}$, there exists a unique $k_i \in \{1,2,...,p\}$ such that 
$$ |g^{\tau}_{k_i} (z_*(t)) - x_i(t)| \leq \epsilon_o \leq \bar{\epsilon}/8 \qquad \forall t \in [\tau,\tau+\delta].  $$

\item For each $i \in \{1,2,...,p\}$, there exists a unique $l_i \in \{1,2,...,p\}$ such that 
$$ |g^{\tau}_{i} (z_*(t)) - x_{l_i}(t)| \leq \epsilon_o \leq \bar{\epsilon}/8 \qquad \forall t \in [\tau,\tau+\delta].  $$
\end{itemize}
A consequence of the latter two items and \eqref{eqsepar} is that 
\begin{equation} \label{eqz1bis}
\hspace{-0.3cm} |x_j(t) - x_{i}(t)| \geq 3\bar{\epsilon} / 4 \quad \forall t \geq \tau_1, ~ \forall i \neq j \in \{1,2,...,p\}.  
\end{equation}
Now, using \eqref{eqaddednew}, we conclude the existence of $\tau_2 \geq \tau_1$ such that, for each $\tau \geq \tau_2$, there exists a unique  $i \in \{1,2,...,p\}$ such that 
$$ | \hat{x}(t) - g^{\tau}_{i} (z(t)) | \leq \epsilon_o  \qquad \forall t \in [\tau,\tau+\delta].  $$
Thus, for each $\tau \geq \tau_2$, there exists a unique  $l \in \{1,2,...,p\}$ such that 
$$ |x_{l}(t) - \hat{x}(t)| 
\leq 2 \epsilon_o \leq \bar{\epsilon}/4 \qquad \forall t \in [\tau,\tau+\delta].  $$
Combining the latter inequality to \eqref{eqz1bis} and using the continuity of solutions and the continuity of the selection $\hat{x}$, we conclude that $l$ is invariant with $\tau$, i.e., there exists a unique  $l \in \{1,2,...,p\}$ such that 
$$ |x_{l}(t) - \hat{x}(t)| \leq 2\epsilon_o \qquad \forall t \in [\tau_2, + \infty).  $$
\end{pf}

\begin{figure*}
    \centering
    \includegraphics[width=\textwidth]{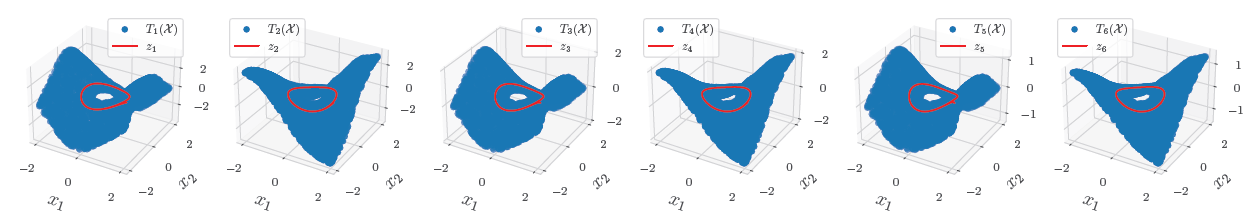}
    \caption{Image sets $T_i(\tilde{\X})$ with $\tilde{\X}=[-2,2]^2\setminus \mathcal{B}(0,0.5)$ and trajectories $t\mapsto z_i(t)$.}
    \label{fig:z_imageT}
\end{figure*}

\begin{example}
In Example \ref{ex:nolcos}, if we pick $\tilde{\X}$ to be a compact subset of $\RR^2\setminus \{0\}$, Assumptions \ref{ass1a} and \ref{ass1b} are satisfied. Besides, if $\tilde{\X}$ contains in its interior the circle $\C$ centered at $0$ and of radius $1$, then any solution $x$ with $x(0) \neq 0$ converges to $\C$ and is eventually in $\tilde{\X}$. Then, Assumption \ref{asssy}  holds since $x$ and $-x$ are the only indistinguishable solutions producing the same output $y$ and not converging to each other. On the other hand, Assumption \ref{asszz} does not hold, as it can be seen in Figure \ref{fig:z_imageT}, since the solutions $t\mapsto z(t)$ eventually circle around the ``hole'' in $T(\tilde{\X})$ obtained by removing a neighborhood of $0$ in $\X$. Nevertheless, Assumption \ref{assz1} holds because the time of revolution around the hole is lower bounded by a positive time $\delta^*$, so that, even in the case of a periodic solution, $z([\tau,\tau + \delta)$ with $\delta<\delta^*$ is always simply connected. The convergence observed in \cite{PaulineNolcos} in thus justified by applying Theorem \ref{thm4}.
\end{example}

\section{Discussion on Assumptions} \label{Sec.VI}

\subsection{Links Between Assumptions \ref{assdist} and \ref{asssy} Under Assumptions \ref{ass1a} and \ref{ass1b}}

\begin{proposition}
Assume system \eqref{sys} is analytic and let $T:\O \subset\RR^{n_x} \to \RR^{n_z}$ be $C^1$. Consider a compact set $\X \subset \O$ such that  Assumptions \ref{ass1a}, \ref{ass1b}, and \ref{assdist} hold. Let  $y:\RR_{\geq 0}\to \RR^{n_y}$ be an output to \eqref{sys} generated by a solution $x_1$ that remains in $\X$ after some $\bar{t} > 0$. 
Then, the following properties are true:
\begin{itemize}
\item There exist at most $p$ solutions $\{x_i\}_{i=1}^p:[\bar{t},
\infty)\to \X$ to \eqref{sys} generating $y$. Besides, such solutions cannot converge to each other. 

\item There are exactly $p$ such solutions provided that one of the following holds: 
\begin{enumerate}
    \item  The set $\X$ is forward invariant.
    \item There exists $\mathcal{U} \subset \mathbb{R}^{n_x}$ a neighborhood of $\X$ such that
    $$ \ind_\mathcal{U}(x_1(t)) = \ind_\X(x_1(t)) \qquad \forall t \geq \bar{t}. $$
    \item The two following properties hold
    \begin{align*}
    \{f(x), -f(x)\} \cap T_{\X} (x) & \neq \emptyset \qquad \forall x \in T^-(x_1(\bar{t})), \\
     \{-f(x), f(x) \} \cap T_{\partial \X}(x) & = \emptyset \qquad \forall x \in \partial \X. 
    \end{align*}
\end{enumerate}
\end{itemize}
\end{proposition}

\begin{pf}
To prove the first item, given $\bar{t} > 0$, we consider two distinct solutions $x_1,x_2:[\bar{t},\infty)\to \X$ such that $$ h(x_1(t))=h(x_2(t))=y(t) \qquad  
\forall  t \in [\bar{t},+\infty). $$ 
By uniqueness of solutions, $x_1(t)\neq x_2(t)$ for all $t\geq \bar{t}$. Also, by analityticity of \eqref{sys}, $x_1(t)\indist x_2(t)$ for all $t\geq \bar{t}$ and  by Assumption \ref{assdist}, $T(x_1(t))=T(x_2(t))$ for all $t\geq \bar{t}$. It follows that for all $t\geq \bar{t}$, $x_1(t)$ and $x_2(t)$ are in $T^-(T(x_1(t))$ and $T^-$ has cardinal $p$ by Assumption \ref{ass1a}. It follows that there are at most $p$ such solutions. Now, let us show that any such solutions cannot converge to each other. Indeed, assume $x_1-x_2$ converges to zero. Since $\X$ is compact, there exists an increasing and diverging sequence of times $(t_n)_{n\in \NN}$  and a limit point $x^*\in \X$ such that $\lim_{n\to +\infty} x_1(t_n)=x^*$. Then, also, $\lim_{n\to +\infty} x_2(t_n)=x^*$. Under Assumption \ref{ass1b}, there exists a neighborhood of $x^*$ where $T$ is locally injective, which contradicts $x_2(t_n)\in T^-(T(x_1(t_n))$ for sufficiently large $n$.

Now, to prove the second item, we assume 
that $y$ is generated, on some interval $[\bar{t},+\infty)$, by a solution $x_1:[\bar{t},+\infty)\to \X$. Under Assumption \ref{ass1a}, we let 
$$ \{x_{1,o},...,x_{p,o}\} := T^-(T(x_1(\bar{t}))) ~ \text{with} ~ x_{1,o}=x_1(\bar{t}). $$ 
 By Assumption \ref{assdist}, we conclude that
$$ x_{i,o} \indist x_{j,o} \qquad  \forall i,j \in \{1,\ldots, p\}. $$ 
Now, for each $i \in \{1,2,...,p \}$, we let the (unique) maximal solution $x_i : I_i \to \RR^{n_x}$ with $\bar{t} \in I_i$ such that $x_i(\bar{t}) =: x_{i,o}$. By definition of $\indist$, for each $i \in \{1,2,...,p\}$, we have 
$$ h(x_i(t))=h(x_1(t)) \qquad \forall t\in I_i\cap I_1 \cap (-\infty, \bar{t}].
$$ 
By analyticity, actually, we conclude that $$ h(x_i(t))=h(x_1(t)) \qquad  \forall t\in I_i\cap I_1. $$ 
Denote $I := \cap_{i\in \{1,\ldots,p\} }I_i$, which is an open set that contains $\bar{t}$. We; thus, have 
$$ x_i(t)\indist x_j(t) \qquad  \forall t\in I, \quad \forall i,j\in \{1,\ldots, p\}. $$ 
To complete the proof, we need to show that 
\begin{align} \label{eqremain}
 x_i(t)\in \X \qquad  \forall t \geq \bar{t}, \quad \forall i\in \{1,\ldots, p\}. 
 \end{align}
Under item 1), when $\X$ is forward invariant, \eqref{eqremain} follows. 

Furthermore, we show that \eqref{eqremain} holds under item 2) using contradiction.
That is, let $t_\X \geq \bar{t}$ be the maximal time, such that, all the $x_i$s remain in $\X$ on $[\bar{t},t_\X]$; namely,
$$ t_\X := \max \{t \in I : x_i(s) \in \X, \ \forall s\in[\bar{t},t], \ \forall i \in \{1,\ldots,p\} \}. $$
Furthermore, let $k \in \{2,3,...,p\}$ be such that the solution $x_k$ leaves the set $\X$ at $t_\X$. Hence, since $x_k$ is continuous, we conclude the existence of $t_m > t_\X$ such that $x_k(t_m) \in \mathcal{U} \backslash \X$. 
By analyticity, we conclude that 
$$ h(x_k(t)) = y(t)  \qquad \forall t \in [\bar{t},t_m].  $$
The latter plus the fact that
$x_k(\bar{t}) \in \ind_{\mathcal{U}}(x_1(\bar{t})) $
 imply that
$x_k(t_m) \in \ind_{\mathcal{U}}(x_1(t_m))$. However, using item 2), we know that  $\ind_{\mathcal{U}}(x_1(t_m)) = \ind_{\mathcal{X}}(x_1(t_m))  \subset \X$, which yields a contradiction.

Suppose, now, that item 3) holds and let  $t_\X$ be finite. By compactness of $\X$, there exists $\varepsilon>0$ such that, for each $i,j\in \{1,\ldots, p\}$ with $i\neq j$, we have 
\begin{equation}
\label{eq_dist_ij_epsilon}
|x_i(t)-x_j(t)| \geq \varepsilon \qquad \forall t\in [\bar{t},t_\X],
\end{equation}
and by Assumption \ref{ass1a}, we conclude that
$$
T^-(T(x_1(t))) = \{x_1(t),\ldots,x_p(t)\} \in \X \quad \forall t\in [\bar{t},t_\X].
$$
Now, by definition of $t_\X$, we conclude the existence of $k \in \{2,3,...,p\}$  such that  $x_k(t_\X) \in \partial \X$ and $x_k$ leaves $\X$ right after $t_\X$. Hence,
using Lemma \ref{LemNag1} we conclude that
$f(x_k(t_\X)) \in T_{\mathbb{R}^{n_x} \backslash \X} (x_k(t_\X))$.
Now, since 
$ f(x_k(t_\X)) \notin T_{\partial \X} (x_k(t_\X))$, we use Lemma \ref{LemNag3} to conclude that
$f(x_k(t_\X)) \in \mathbb{R}^{n_x} \backslash  T_{\X} (x_k(t_\X))$. Hence, using Lemma \ref{LemNag2}, we conclude that, for some $\delta > 0$ sufficiently small, 
\begin{align} \label{eqleave} 
x_k(t) \notin \X \qquad \forall t \in (t_\X, t_\X+\delta] \subset I. 
\end{align}
Note that this $k$ may not be  unique; however, for simplicity we let it be unique in this proof, the exact same arguments apply in the general case. Now, we distinguish between two scenarios:

\begin{itemize} 
\item When $t_\X > \bar{t}$, we conclude that $x_k([\bar{t}, t_\X]) \subset \X$. Hence, using Lemma \ref{LemNag1}, we conclude that 
$$-f(x_k(t_\X)) \in T_{\X}(x_k(t_\X)). $$
However, since 
$-f(x_k(t_\X)) \notin T_{\partial \X}(x_k(t_\X))$, we use Lemma \ref{LemNag3} to conclude that  $-f(x_k(t_\X)) \in 
D_{\X}(x_k(t_\X))$. Hence, using Lemma \ref{LemNag2}, we conclude that for $\delta > 0$ sufficiently small, we have  
\begin{align} \label{eqbackenter}
x_k(t) \in \text{int}(\X) \qquad  
\forall t\in [t_\X-\delta, t_\X)
\subset I. 
\end{align}

\item When $t_\X = \bar{t}$, we use \eqref{eqleave}, to conclude that  
$$ f(x_k(\bar{t})) \in T_{\mathbb{R}^{n_x} \backslash \X}(x_k(\bar{t})). $$  
Furthermore, since $ f(x_k(\bar{t})) \notin T_{\partial \X}(x_k(\bar{t})) $; we conclude, using Lemma \ref{LemNag3} that 
$ f(x_k(\bar{t})) \notin T_{ \X}(x_k(\bar{t}))$. As a result, according to item 3), we have  
$$ - f(x_k(\bar{t})) \in T_{\X}(x_k(\bar{t})) \backslash T_{\partial \X}(x_k(\bar{t})) = D_\X(x_k(\bar{t})). $$
Hence, using Lemma \ref{LemNag2}, \eqref{eqbackenter} follows. 
\end{itemize}

Now, using continuity of $T^-$ and Assumption \ref{ass1a}, we conclude that, for each $t\in (t_\X, t_\X+\delta]$, there exists $x_{k,t}'\in \X$ such that
$$
T^-(T(x_1(t)) = \{ x_1(t),\ldots,x_{k-1}(t), x_{k,t}', x_{k+1}(t),\ldots,x_{p}(t) \},
$$
and
$$ \lim_{t\to t_\X^+} x_{k,t}' = x_k(t_\X). $$
Furthermore, for $\delta > 0$ even smaller, we conclude that all the elements of $T^-(T(x_1(t))$, for all $t \in [t_\X - \delta, t_\X + \delta]$, are at a distance larger than $\varepsilon/2$ from each other. 
By continuity of solutions with respect to initial data, we conclude the existence $t_m\in (t_\X, t_\X+\delta]$ sufficiently close to $t_\X$ 
 such that the solution $x_k'$ starting from $x_{k_m,t_m}'$ at time $t_m$ verifies 
 $$ \|x_k(t)-x_k'(t)\|< \varepsilon/2 \qquad \forall  [t_\X-\delta,t_\X+\delta] $$ 
 and $ x_k'(t_\X-\delta) \in \text{int}(\X)$.
 But by Assumption \ref{assdist}, $x_{k,t_m}'\indist x_1(t_m)$ and therefore, by analyticity, we have 
 $$ x_k(t_\X-\delta)\indist x_k'(t_\X-\delta). $$ 
 It follows from \eqref{eq_dist_ij_epsilon} that $x_k(t_\X-\delta)=x_k'(t_\X-\delta)$ and therefore, $x_k=x_k'$ by uniqueness of solutions. This is impossible since $x_k(t_m)\notin \X$ while $x_k'(t_m)=x_{k,t_m}'\in \X$.
\end{pf}

The previous proposition shows that, when system \eqref{sys} is analytic and under some assumptions on the way solutions might exit or enter $\X$, Assumptions \ref{ass1a}, \ref{ass1b}, and \ref{assdist} 
imply Assumption \ref{asssy}.  
In the following result, we investigate a converse statement. That is, for $T$ solution to \eqref{eq_PDE} satisfying Assumptions \ref{ass1a} and \ref{ass1b}, we consider an output $y$ satisfying the following assumption, which is slightly stronger than Assumption \ref{asssy}.

\begin{assumption} \label{asssybis}
Given $y:\RR_{\geq 0}\to \RR^{n_y}$, there exist $\bar{t}>0$ and at least $p \in \mathbb{N}_{>0}$ solutions $ \{ x_i \}^p_{i=1}: [\bar{t},+\infty) \to \X$ to \eqref{sys} such that for all $i\in \{1,\ldots,p\}$, 
\begin{itemize}
    \item For each $t\geq \bar{t}$, $h(x_i(t))=y(t)$. 
    \item For each $i$, $j \in \{1,\ldots,p\}$, we have 
\begin{align} \label{eqliminf}    
\liminf_{t \rightarrow \infty} |x_i(t)-x_{j\neq i}(t)| > 0.
\end{align}
\end{itemize}
\end{assumption}
As a result, we show that $T^-(T(x)) \subset  \ind_\X(x)$, for all $x$ within the $\omega$-limit set of the solutions in Assumption \ref{asssybis} generating $y$. 

\begin{proposition}
Consider a $C^1$ map $T: \O \subset\RR^{n_x} \to \RR^{n_z}$ verifying \eqref{eq_PDE} and Assumptions \ref{ass1a} and \ref{ass1b} on a compact set $\X \subset \O$. Let $y:\RR_{\geq 0} \to \RR^{n_y}$ be an output to \eqref{sys} such that Assumption \ref{asssybis} holds.
Let $x_1$ be one of the solutions introduced in Assumption \ref{asssybis} that generates $y$. Then, 
$$ T^-(T(x)) \subseteq \ind_{\X}(x) \qquad \forall x \in \omega(x_{1}(\bar{t})). $$
\end{proposition}

\begin{pf}
Let the solutions $\{x_i\}_{i=1}^p : [\bar{t}\to \infty)\to \X$ to \eqref{sys} generating $y$. Then, $x_{\rm ext}:=(x_1,x_2,\ldots,x_p)$ is solution to a duplicated system with dynamics $f_{\rm ext}=(f,f,\ldots,f)$ and verifies 
\begin{equation}
\label{eq_hequals}
    h(x_i(t))=h(x_j(t)) \quad \forall t\geq \bar{t} ~ \forall (i,j)\in 
    \{1,\ldots,p\}^2.
\end{equation}
Let $\Omega := \omega(x_{\rm ext}(\bar{t}))$,
its $\omega$-limit set.
Since $\X$ is compact and by uniqueness of solutions, $\Omega$ is compact, contained in $\X$, and invariant by $f_{\rm ext}$. 
Since $h$ is continuous, we deduce from \eqref{eq_hequals} that for all $(x_1^\star,\ldots,x_p^\star)\in \Omega$,
\begin{equation}
\label{eq_hequal_star}
    h(x_i^\star)=h(x_j^\star) \qquad \forall (i,j)\in \{1,\ldots,p\}^2 .
\end{equation}
Now pick $x_{\rm ext}^\star:=(x_1^\star,\ldots,x_p^\star)\in \Omega$.
By backward-invariance of $\Omega$, the backward solution $(\phi_1, ..., \phi_p)$ by $f_{\rm ext}$ initialized at $x_{\rm ext}^\star$ is defined on $\RR_{\leq 0}$ and remains in $\Omega$. It follows from \eqref{eq_hequal_star} (which holds everywhere on $\Omega$) that 
$h(\phi_i(t))=h(\phi_j(t))$ for all 
$t \leq 0$ and for all 
$ (i,j) \in \{1,\ldots,p\}^2$.

Since for each $i\in \{1,\ldots,p\}$, $\phi_i$ is the backward solution to \eqref{sys} initialized at $x_i^\star$ by definition of $f_{\rm ext}$, we deduce that
$ x_1^\star \indist x_2^\star \indist ... \indist x_p^\star$, and thus, since each $x_i^\star\in \X$, $x_i^\star\in \ind_\X(x_1^\star)$.
Moreover, by uniqueness of the solutions to \eqref{eq_PDE} on backward-invariant compact sets \cite[Theorem 2]{BriAndBerSer}, we know that, for each $i \in \{1,2,...,p \}$, 
$ T(x_i^\star) = \int_{-\infty}^0 e^{-As} B h(\phi_i(s))ds, $
which implies that 
$ T(x_1^\star)= ... = T(x_p^\star) $.
The latter allows us to conclude that
\begin{align} \label{eqequal-}
\{x_1^\star,\ldots,x_p^\star\} \subseteq T^-(T(x_1^\star)) \cap \ind_\X(x_1^\star). 
\end{align}
Now, under the second item in Assumption \ref{asssybis},   the $x_i^\star$s are distinct, then necessarily, we would have 
\begin{align} \label{eqequal}
\{x_1^\star,\ldots,x_p^\star\} = T^-(T(x_1^\star)) \subseteq \ind_\X(y_1). 
\end{align}
\end{pf}

When, instead of \eqref{eqliminf}, we assume that the solutions $\{x_i\}^p_{i=1}$ in Assumption \ref{asssybis} do not converge to each other (as in Assumption \ref{asssy}),
we cannot conclude the equality in \eqref{eqequal} based on \eqref{eqequal-}. 
Indeed, it is possible to have, for example, the solutions $x_1$ and $x_2$ converge to a same $\omega$-limit point $x_1^* = x_2^*$ along a same sequence of times, and they still do not converge to each other; see Example \ref{exp44}. 

\begin{example} \label{exp44}
Consider the two-dimensional 
system 
$$ 
\dot{x} = f(x) := \begin{bmatrix} 
(1 -|x|^2) x_1 + x_2 
\\
(1 - |x|^2) x_2 - x_1
\end{bmatrix} ~~ x \in \X := \{ |x| \leq 2 \}.
$$ 
Note that this system admits a unique nontrivial limit cycle describing the unite circle 
$\{ |x| = 1 \}$. The latter  attracts all the solutions except the one starting from the origin $\{x = 0\}$.
Next, we propose to re-scale the system using the smooth function $\phi : \mathbb{R}_{\geq 0} \rightarrow [0,1]$ given by
$$ \phi(r) := 
\left\{ 
\begin{matrix} 
0 &  \text{if} ~ r \geq 1 
\\
(1 - r)^2 & \text{if} ~ r < 1.
\end{matrix}
\right.     $$
As a result, we introduce the system  
\begin{align} \label{sysphif}
\dot{x} = \phi(|x|^2) f(x) \qquad x \in \X,
\end{align}
for which, every point within the set $\{ |x| = 0 \} \cup \{ |x| \geq 1 \}$ is a static equilibrium point, and the solutions starting from the set $\{0 < |x| < 1\}$ converge to the set  $\{|x| \geq 1\}$ by spiraling, i.e., without converging to any specific point in $\{|x| \geq 1\}$. Hence, by letting $x_1$ be a solution to \eqref{sysphif} starting from $\{ |x| =1 \}$ (which remains there, i.e., $x_1(t) = x_1(0)$ for all $t \geq 0$) and $x_2$ be a solution to \eqref{sysphif} starting from $ \{0 < |x| < 1\}$, we conclude that they do not converge to each other (since $x_2$ keeps spiralling), but still they share a common $\omega$-limit point which is $x_1(0)$.    
\end{example}

\subsection{Link Between Assumption \ref{ass1b} and the Rank of a Differential Observability Map}
\label{sec_diff_obs}

In this section, we assume $f$ and $h$ to be smooth. Given a positive integer $m$, consider the map $H_m$ defined by
\begin{equation} \label{eq_def_H}
H_m(x) := \left( h(x), L_fh(x),\hdots, L_f^{m-1}h(x) \right),
\end{equation}
containing the output map $h$ and its $m-1$ Lie derivatives.
The injectivity of $H_m$ characterizes the so-called \textit{differential observability} of order $m$, meaning that the state is uniquely determined from the knowledge of the output and its first $m-1$ time derivatives. In this paper where we handle non-observable systems, we will not make such an assumption. However, the next result shows that the rank of the KKL map $T$ may be related to the rank of the map $H_m$, when the eigenvalues of $A$ are picked sufficiently fast. This result allows to characterize the set of regular points where the previous convergence results hold by simply checking the rank of $H_m$. 

\begin{proposition}
\label{prop_diff_obs}
Consider a controllable pair $(A_o,B_o)\in \RR^{m\times m}\times \RR^{m}$ with $A_o$ Hurwitz and $m\in \NN$ such that the map $H_m$ defined in \eqref{eq_def_H} is full-rank on the compact set $\X$. Define $
A = I_{n_y} \otimes A_o,
B = I_{n_y} \otimes B_o$.
There exists $k^\star>0$ such that for all $k>k^\star$, there exists a $C^1$ map $T_k:\RR^{n_x}\to \RR^{n_z}$ which is full-rank on $\X$ and verifies \eqref{eq_PDE} with the pair $(kA,B)$.
\end{proposition}

\begin{pf}
First, since in \eqref{eq_PDE}, the map $f$ is only used on $\X$, we can consider an open bounded set $\O$ and a modified vector field $\breve{f}$ such that $f= \breve{f}$ on $\X$ and $\O$ is backward invariant for $\dot{x} = \breve{f}(x)$. Then, as in Theorem \ref{thm_info}, we show that the map defined by
$$
T_k(x) := \int_{-\infty}^0 e^{-kAs} B h(\breve{X}(x,s)) ds, 
$$
with $\breve{X}(x,s)$ representing the flow operator for $\dot{x} = \breve{f}(x)$, is $C^1$ on $\O$ for $k$ sufficiently large and verifies \eqref{eq_PDE} for the pair $(kA,B)$. Given the structure of $A$ and $B$, $T_k(x)= (T_{k,1}(x),\ldots,T_{k,n_y}(x))$, where for all $i\in \{1,\ldots, n_y\}$, 
$$
T_{k,i}(x) = \int_{-\infty}^0 e^{-kA_os} B_o h_i(\breve{X}(x,s)) ds.
$$
Then, after $m$ integration by parts, we see that
$$
T_{k,i}(x) = A_o^{-m}\C_o K\left( H_{m,i}(x) + \frac{1}{k^{m}} K^{-1}\C_o^{-1} R_{k,i}(x) \right)
$$
where $\C_o$ is a square controllability matrix associated to $(A_o,B_o)$ (thus invertible), $K=\diag\left(\frac{1}{k}, \ldots, \frac{1}{k^{m}} \right)$, $H_{m,i}$ is the observability map of order $m$ associated to $h_i$ and $R_{k,i}:\RR^{n_x}\to \RR^{m}$ are the remainders defined as
$$
R_{k,i}(x) := \int_{-\infty}^0 e^{-kA_os} B_o L_f ^{m} h_i(\breve{X}(x,s)) ds.
$$
Defining $R_k(x)= (R_{k,1}(x),\ldots,R_{k,n_y}(x))$, and using similar arguments as in proving $T_k$ is $C^1$, for $k$ sufficiently large, $R_k$ is $C^1$ and
\begin{equation*}
  \frac{\partial T_k}{\partial x}(x) = M_1(k) \left(\frac{\partial H_{m}}{\partial x}(x) + M_2(k)\frac{\partial R_k}{\partial x} (x) \right),
\end{equation*}
where $M_1(k) :=[I_{n_y}\otimes (A_o^{-m}\C_o K)]M$, $M_2(k) :=\frac{1}{k^{m}} M^{-1}[I_{n_y}\otimes (K^{-1}\C_o^{-1})]$, and $M$ is an invertible square matrix reordering the lines.
Since $H_{m}$ is full-rank on $\X$ compact, there exists $\rho_H >0$ such that for all $x\in \X$,
$$
\frac{\partial H_{m}}{\partial x}(x)^\top \frac{\partial H_{m}}{\partial x}(x) \geq \rho_H I_{n_z}.
$$
By invertibility of $A_o^{-m}\C_o$ and $M$, and given the structure of $k$, there also exists $\rho_c>0$ such that $M_1(k)^\top M_1(k) \geq \frac{\rho_c}{k^m} I_{n_z}$.
It follows that for all $v\in \RR^n$ and all $x\in \X$, $\left|\frac{\partial H_{m}}{\partial x}(x) v \right|\geq \rho_H |v|$, and 
$$
\left|\frac{\partial T_k}{\partial x}(x) v \right|\geq \frac{\rho_c}{k^{m}}\left(\rho_H |v| -  \left|M_2(k) \frac{\partial R_k}{\partial x} (x) v\right|\right)
$$
Upper-bounding $M_2(k)$ by $k^{m} c_c$ for some $c_c>0$ independent from $k$, we thus get
$$
\left|\frac{\partial T_k}{\partial x}(x) v \right|\geq \frac{\rho_c}{k^{m}}\left(\rho_H |v| - c_c \left|\frac{\partial R_k}{\partial x} (x) v\right|\right)
$$
Let us now upper-bound the jacobian of $R_k$,
$$
\left|\frac{\partial R_{k}}{\partial x}(x) \right | \leq  \int_{-\infty}^0{\left|e^{-kAs}\right| \left|B \frac{\partial L_f^{m}h}{\partial x}(\breve{X}(x,s)) \right| \left|\frac{\partial \breve{X}}{\partial x}(x,s)\right|} .
$$
$A$ being Hurwitz, there exists a positive definite  matrix $P$ and a positive scalar $a$ such that 
$
A^T P + P A \leq -2aP.
$
It follows that for all $s$ in $(-\infty,0]$
$
\left|e^{-kAs}\right|\leq \sqrt{\frac{\lambda_{\max}(P)}{\lambda_{\min}(P)}} e^{k as},
$
where $\lambda_{\min}$ and $\lambda_{\max}$ denote the minimal and maximal eigenvalues respectively.
The bounded set $\O$ being backward invariant under the modified dynamics, we can consider 
$
c_0 = \max_{x\in \cl(\O)} \left|B \frac{\partial L_f^{n_z} h}{\partial x}(x) \right|
$.
Finally, we show that the Jacobian of the flow operator $\frac{\partial \breve{X}}{\partial x}(x,s)$ grows at most exponentially in time. To do so, we use the fact that $\psi(t)=\frac{\partial \breve{X}}{\partial x}(x,t)$ satisfies the ODE :
$
\dot{\psi}(t) = \frac{\partial f}{\partial x}(\breve{X}(x,t)) \psi(t).
$
Denoting 
$c_1 := \max_{x\in\cl(\O)} \left|\frac{\partial \breve{f}}{\partial x}(x)\right|
$
and bounding the previous differential inequality, we obtain that $\left|\frac{\partial \breve{X}}{\partial x}(x,s)\right| \leq e^{-c_1 s}$ for all $s$ in $(-\infty,0]$. Therefore, for $k> \frac{c_1}{a}$, there exists a positive constant $c_2$ such that 
$
\left|\frac{\partial R_k}{\partial x}(x) \right |\leq\frac{c_2}{ka-c_1} 
$. It follows that, for all $(x,v) \in \X \times \RR^n$,
$$
\left|\frac{\partial T_k}{\partial x}(x) v \right|\geq \frac{\rho_c}{k^{n_z}}\left(\rho_H  - c_c \frac{c_2}{ka-c_1} \right)|v|, $$
and $T_k$ is indeed full-rank on $\X$ for $k$ sufficiently large.
\end{pf}

\begin{remark}
Note that if we add the backward invariance condition used in Theorem \ref{thm_info} and we ask from the PDE in \eqref{eq_PDE} to hold on $\O$, the full-rank map $T_k$ given by Proposition \ref{prop_diff_obs} is the unique solution to \eqref{eq_PDE} (see Remark \ref{rem_uniqueness}). However, we cannot directly prove that $T_k$ characterizes the backward distinguishable points, i.e., verifies \eqref{eq_injectivity_indistinguishable}. According to Theorem \ref{thm_info} when $m\leq n_o$, \eqref{eq_injectivity_indistinguishable} holds if $(kA_o,B_o)$ is outside a zero-measure set.
\end{remark}

\begin{example}
For the system in Example \ref{ex:nolcos}, we can see that $H_1(x) := (x_1^2-x_2^2,x_1x_2)$ is full-rank for all $x \neq (0,0)$. On the other hand, $H_m$ is never full-rank at $(0,0)$ for any $m \in \NN$. This means that, at least for $n_o\geq 1$ and $(A_o,B_o)\in \RR^{n_o\times n_o}\times \RR^{n_o}$, the solution to \eqref{eq_PDE} with pair $(k I_{n_y}\otimes A_o,I_{n_y}\otimes B_o)$ should be full-rank on any compact set $\tilde{\X} \subset \RR^2 \setminus 
\{(0,0)\}$, provided that $k$ is sufficiently large. Nothing can be said for compact sets that include $(0,0)$.
\end{example}

\section{Conclusion}
This paper proposed a set-valued KKL observer for a nonlinear system \eqref{sys} disobeying the backward-distinguishability property. Provided that the transformation $T$ is full rank on $\X$  and its preimage has a constant cardinality, it is shown that the set-valued preimage is locally Lipschitz and admits a locally Lipschitz extension $\Tinv$. When, additionally, a given output $y$ is generated by a number of distinct solutions (that equals the cardinality of the preimage) not converging to each other, the designed observer is shown to asymptotically reconstructs each of such solutions. Such an approach applies to any nonlinear system, with no particular normal form, unlike high-gain or sliding-mode based methods in \cite{MorMujEsp,MorBes}.
In the future, it would be interesting to relax some of the considered assumptions and to test this approach on practical examples.   

\appendix 

\section{Some useful lemmas}

\begin{lemma} [Constant-rank theorem] \label{lem1}
Consider a subset $\mathcal{X} \subset \mathbb{R}^{n_x}$, $x_o \in \mathcal{X}$ and a continuously differentiable map $T : \mathcal{X} \rightarrow \mathbb{R}^{n_z}$ such that 
$ \text{rank} \left( \frac{\partial T}{\partial x}(x) \right) = n_x$ for all $x \in \mathcal{X}$.
Then,  there exist a neighborhood of $x_o$ denoted $U(x_o)$, $\phi : U(x_o) \rightarrow \phi(U(x_o))$, a neighborhood of $T(x_o)$, denoted $U(T(x_o))$, and 
$\psi : U(T(x_o)) \rightarrow \psi(U(T(x_o)))$ continuously differentiable diffeomorphisms such that 
$$ \psi \circ T (x)  = (\phi(x), 0,0,...,0)  \qquad \forall x \in U(x_o).  $$
\end{lemma}

Next, we recall the following extension theorem for Lipschitz maps \cite[Theorem 1.1]{Gob} and \cite[Theorem 1.7]{De_Lellis_2011}.

\begin{lemma} [Lipschitz extension theorem] \label{lem2}
Let $\Z \subset \mathbb{R}^{n_z}$
 and $ \tilde{F} : \Z \rightarrow \A_p(\mathbb{R}^{n_x})$  be  Lipschitz.  Then,  there exists a Lipschitz extension 
$ \tilde{F}_a : \mathbb{R}^{n_z} \rightarrow \A_p( \mathbb{R}^{n_x})$ of $\tilde{F}$.  
\end{lemma}

Next, we recall from \cite{Aubin:1991:VT:120830}  some useful results on forward invariance of a set $\X$. The next result can be found in \cite[Proposition 3.4.1]{Aubin:1991:VT:120830}. 

\begin{lemma} \label{LemNag1}
Consider system \eqref{sys}
with $f: \O \rightrightarrows \mathbb{R}^n$ continuous.  Consider a closed set $\X \subset \O$ and a solution $x$ of $\dot{x}\in f(x)$ satisfying 
$ \forall \bar{t} > 0, ~ \exists t \in (0, \bar{t}] : x(t) \in \X. $  
Then,  $f(x(0)) \in T_\X(x(0))$.
\end{lemma}

The next result can be found  in \cite[Theorem 4.3.4]{Aubin:1991:VT:120830}. 

\begin{lemma}  \label{LemNag2}
Consider system \eqref{sys}
with $f: \O \rightrightarrows \mathbb{R}^n$ continuous.  Let $\X \subset \O$ be closed with nonempty interior and $x_o \in \partial \X$. If $f(x_o) \subset D_{\X}(x_o)$, then, for each solution $x$ to $\dot{x}\in f(x)$ starting from $x_o$, there exists $\delta > 0$ such that
$x((0,\delta]) \subset \mbox{int}(\X)$.     
\end{lemma} 

The following result can be found in \cite[Lemma 4.3.2 and Theorem 4.3.3]{Aubin:1991:VT:120830}.

\begin{lemma} \label{LemNag3}
Given a closed set $\X \subset \mathbb{R}^{n_x}$, for each $x \in \partial \X$, we have 
$$ D_\X (x) = \mathbb{R}^{n_x} \backslash T_{\mathbb{R}^{n_x} \backslash \X}(x), ~~
T_{\partial \X} (x) =  T_{\X}(x) \cap T_{\mathbb{R}^{n_x} \backslash \X}(x). $$

\end{lemma}

\bibliographystyle{plain}
\bibliography{biblio_journal}

\begin{thebibliography}{10}

\bibitem{PaulineNolcos}
V.~Alleaume and P.~Bernard.
\newblock {KKL} set-valued observers for non-observable systems.
\newblock {\em IFAC Symposium on Nonlinear Control Systems}, 2022.

\bibitem{Alm}
F.~J. Almgren.
\newblock {\em Almgren’s big regularity paper. Q-Valued Functions Minimizing
  Dirichlet's Integral and the Regularity of Area-Minimizing Rectifiable
  Currents Up to Codimension 2}, volume~1.
\newblock World Scientific Monograph Series in Mathematics, World Scientific
  Publishing Co. Inc., River Edge, NJ, 2000.

\bibitem{andrieuExistenceKazantzisKravaris2006a}
V.~Andrieu and L.~Praly.
\newblock On the {{Existence}} of a {{Kazantzis--Kravaris}}/{{Luenberger
  Observer}}.
\newblock {\em SIAM Journal on Control and Optimization}, 45(2):432--456, 2006.

\bibitem{Aubin:1991:VT:120830}
J.~P. Aubin.
\newblock {\em Viability Theory}.
\newblock Birkhauser Boston Inc., Cambridge, MA, USA, 1991.

\bibitem{aubinDifferentialInclusionsSetValued1984a}
J.~P. Aubin and A.~Cellina.
\newblock {\em Differential {{Inclusions}}: {{Set-Valued Maps}} and {{Viability
  Theory}}}, volume 264 of {\em Grundlehren Der Mathematischen
  {{Wissenschaften}}}.
\newblock {Springer Berlin Heidelberg}, 1984.

\bibitem{Banach1934}
S.~Banach and S.~Mazur.
\newblock Über mehrdeutige stetige abbildungen.
\newblock {\em Studia Mathematica}, 5(1):174--178, 1934.

\bibitem{bernardObserverDesignContinuoustime2022}
P.~Bernard, V.~Andrieu, and D.~Astolfi.
\newblock Observer design for continuous-time dynamical systems.
\newblock {\em Annual Reviews in Control}, 2022.

\bibitem{BerPra21}
P.~Bernard and L.~Praly.
\newblock Estimation of position and resistance of a sensorless {PMSM} : a
  nonlinear luenberger approach for a non-observable system.
\newblock {\em IEEE Trans. on Automatic Control}, 66:481--496, 2021.

\bibitem{BriAndBerSer}
L.~Brivadis, V.~Andrieu, P.~Bernard, and U.~Serres.
\newblock Further remarks on {KKL} observers.
\newblock {\em To appear in Systems and Control Letters, Available online at
  \url{https://hal.archives-ouvertes.fr/hal-03695863}}, 2022.

\bibitem{nvalued2018}
R.~F. Brown and D.~L. Goncalves.
\newblock On the topology of n-valued maps.
\newblock {\em Adv. Fixed Point Theory}, 8:205--220, 2018.

\bibitem{KKLtoolbox}
M.~Buisson-Fenet, L.~Bahr, and F.~Di-Meglio.
\newblock Learning to observe : neural network-based {KKL} observers.
\newblock {Python toolbox available at}
  \url{https://github.com/Centre-automatique-et-systemes/learn_observe_KKL.git},
  2022.

\bibitem{Gob}
J.~Goblet.
\newblock Lipschitz extension of multiple banach-valued functions in the sense
  of almgren.
\newblock {\em Houston Journal of Mathematics}, 35(1):223--231, 2009.

\bibitem{kazantzisNonlinearObserverDesign1998}
N.~Kazantzis and C.~Kravaris.
\newblock Nonlinear observer design using {{Lyapunov}}’s auxiliary theorem.
\newblock {\em Systems \& Control Letters}, 34(5):241--247, 1998.

\bibitem{KreisselmeierEngel_TAC_03}
G.~Kreisselmeier and R.~Engel.
\newblock {Nonlinear observers for autonomous Lipschitz continuous systems}.
\newblock {\em IEEE Trans. on Automatic Control}, 48(3), 2003.

\bibitem{krener2001nonlinear}
A.J. Krener and M-Q. Xiao.
\newblock Nonlinear observer design in the siegel domain through coordinate
  changes.
\newblock {\em IFAC Proceedings Volumes}, 34(6):519--524, 2001.

\bibitem{De_Lellis_2011}
C.~De Lellis and E.~Spadaro.
\newblock $q$-valued functions revisited.
\newblock {\em Memoirs of the American Mathematical Society}, 211(991), 2011.

\bibitem{luenbergerObservingStateLinear1964}
D.~G. Luenberger.
\newblock Observing the {{State}} of a {{Linear System}}.
\newblock {\em IEEE Trans. on Military Electronics}, 8(2):74--80, 1964.

\bibitem{MorMujEsp}
J.~A. Moreno, H.~Mujica-Ortega, and G.~Espinosa-P\'erez.
\newblock A global bivalued-observer for the sensorless induction motor.
\newblock {\em IFAC World Congress}, 2017.

\bibitem{MorBes}
J.A. Moreno and G.~Besan\c{c}on.
\newblock Multivalued finite-time observers for a class of nonlinear systems.
\newblock {\em IEEE Conference on Decision and Control}, 2017.

\bibitem{Shoshitaishvili_TSA_90}
A.N. Shoshitaishvili.
\newblock {Singularities for projections of integral manifolds with
  applications to control and observation problems}.
\newblock {\em Theory of singularities and its applications}, 1:295, 1990.

\bibitem{staecker2021partitions}
C.~P. Staecker.
\newblock Partitions of $ n $-valued maps.
\newblock {\em arXiv preprint arXiv:2101.09326}, 2021.

\bibitem{VerLorForMenZar}
C.M. Verrelli, E.~Lorenzani, R.~Fornari, M.~Mengoni, and L.~Zarri.
\newblock Steady-state speed sensor fault detection in induction motors with
  uncertain parameters: A matter of algebraic equations.
\newblock {\em Control Engineering Practice}, 80:125--137, 2018.

\end{thebibliography}

\end{document}